\documentclass[3p]{elsarticle}
\usepackage{bm}
\usepackage[parfill]{parskip}
\usepackage[utf8]{inputenc}
\usepackage{url}
\usepackage{amsmath,amssymb,amsfonts,amsthm}
\usepackage{color}
\usepackage{xcolor}
\usepackage{bm}
\usepackage{graphicx}
\usepackage{subcaption}
\usepackage{lipsum}
\usepackage{amsfonts}
\usepackage{graphicx}
\usepackage{epstopdf}
\usepackage{array}
\usepackage[hidelinks]{hyperref}
\usepackage{lineno}
\usepackage{bbm}
\usepackage{xfrac}
\usepackage{varwidth}
\usepackage[nice]{nicefrac}
\usepackage{yhmath}

\usepackage{amsmath}
\usepackage{algorithm}
\usepackage[noend]{algpseudocode}

\ifpdf
  \DeclareGraphicsExtensions{.eps,.pdf,.png,.jpg}
\else
  \DeclareGraphicsExtensions{.eps}
\fi

\newtheorem{theorem}{Theorem}
\newtheorem{lemma}{Lemma}

\newtheorem{definition}{Definition}

\makeatletter
\newcommand{\doublehat}[1]{%
\begingroup%
  \let\macc@kerna\z@%
  \let\macc@kernb\z@%
  \let\macc@nucleus\@empty%
  \widehat{\raisebox{.35ex}{\vphantom{\ensuremath{#1}}}\smash{\widehat{#1}}}%
\endgroup%
}
\makeatother

\makeatletter
\newsavebox\myboxA
\newsavebox\myboxB
\newlength\mylenA

\newcommand*\xoverline[2][0.75]{%
    \sbox{\myboxA}{$\m@th#2$}%
    \setbox\myboxB\null
    \ht\myboxB=\ht\myboxA%
    \dp\myboxB=\dp\myboxA%
    \wd\myboxB=#1\wd\myboxA
    \sbox\myboxB{$\m@th\overline{\copy\myboxB}$}
    \setlength\mylenA{\the\wd\myboxA}
    \addtolength\mylenA{-\the\wd\myboxB}%
    \ifdim\wd\myboxB<\wd\myboxA%
       \rlap{\hskip 0.5\mylenA\usebox\myboxB}{\usebox\myboxA}%
    \else
        \hskip -0.5\mylenA\rlap{\usebox\myboxA}{\hskip 0.5\mylenA\usebox\myboxB}%
    \fi}
\makeatother

\usepackage{tikz}
\usetikzlibrary{chains, positioning, arrows.meta, bending, shapes.arrows}
\usetikzlibrary{positioning,shapes,arrows,calc,shapes.geometric,backgrounds,fit}
\usepackage{pifont}
\usetikzlibrary{positioning,shapes,arrows,calc,shapes.geometric,backgrounds,fit}

\tikzstyle{block} = [rectangle,draw,minimum width=2em,align=center,rounded corners, minimum height=2em,scale=1.0]
\tikzstyle{blockleft} = [rectangle,draw,minimum width=2em,align=left,rounded corners, minimum height=2em,scale=1.0]
\tikzstyle{bigblock} = [rectangle,draw,minimum width=8em,align=center,rounded corners, minimum height=4em,scale=1.0]
\tikzstyle{connect} = [draw,-latex']
\tikzstyle{decision} = [diamond, draw, 
    text width=4.5em, text badly centered, node distance=3cm, inner sep=0pt]
\tikzstyle{line} = [draw, -latex']
\tikzstyle{cloud} = [draw, ellipse,fill=red!20, node distance=3cm,
    minimum height=2em]
\tikzstyle{linenoarrow}=[draw]

\makeatletter
\let\save@mathaccent\mathaccent
\newcommand*\if@single[3]{%
  \setbox0\hbox{${\mathaccent"0362{#1}}^H$}%
  \setbox2\hbox{${\mathaccent"0362{\kern0pt#1}}^H$}%
  \ifdim\ht0=\ht2 #3\else #2\fi
  }
\newcommand*\rel@kern[1]{\kern#1\dimexpr\macc@kerna}
\newcommand*\widebar[1]{\@ifnextchar^{{\wide@bar{#1}{0}}}{\wide@bar{#1}{1}}}
\newcommand*\wide@bar[2]{\if@single{#1}{\wide@bar@{#1}{#2}{1}}{\wide@bar@{#1}{#2}{2}}}
\newcommand*\wide@bar@[3]{%
  \begingroup
  \def\mathaccent##1##2{%
    \let\mathaccent\save@mathaccent
    \if#32 \let\macc@nucleus\first@char \fi
    \setbox\z@\hbox{$\macc@style{\macc@nucleus}_{}$}%
    \setbox\tw@\hbox{$\macc@style{\macc@nucleus}{}_{}$}%
    \dimen@\wd\tw@
    \advance\dimen@-\wd\z@
    \divide\dimen@ 3
    \@tempdima\wd\tw@
    \advance\@tempdima-\scriptspace
    \divide\@tempdima 10
    \advance\dimen@-\@tempdima
    \ifdim\dimen@>\z@ \dimen@0pt\fi
    \rel@kern{0.6}\kern-\dimen@
    \if#31
      \overline{\rel@kern{-0.6}\kern\dimen@\macc@nucleus\rel@kern{0.4}\kern\dimen@}%
      \advance\dimen@0.4\dimexpr\macc@kerna
      \let\final@kern#2%
      \ifdim\dimen@<\z@ \let\final@kern1\fi
      \if\final@kern1 \kern-\dimen@\fi
    \else
      \overline{\rel@kern{-0.6}\kern\dimen@#1}%
    \fi
  }%
  \macc@depth\@ne
  \let\math@bgroup\@empty \let\math@egroup\macc@set@skewchar
  \mathsurround\z@ \frozen@everymath{\mathgroup\macc@group\relax}%
  \macc@set@skewchar\relax
  \let\mathaccentV\macc@nested@a
  \if#31
    \macc@nested@a\relax111{#1}%
  \else
    \def\gobble@till@marker##1\endmarker{}%
    \futurelet\first@char\gobble@till@marker#1\endmarker
    \ifcat\noexpand\first@char A\else
      \def\first@char{}%
    \fi
    \macc@nested@a\relax111{\first@char}%
  \fi
  \endgroup
}
\makeatother

\journal{arXiv.org}

\newcommand{\TheTitle}{An energy stable and conservative multiplicative dynamical low-rank discretization for the Su-Olson problem} 

\date{\today}

\usepackage{amsopn}
\DeclareMathOperator{\diag}{diag}

\journal{arXiv}

\begin{document}
\begin{frontmatter}

\title{\TheTitle}

\author[adressWuerzburg]{Lena Baumann}
\author[adressInnsbruck]{Lukas Einkemmer}
\author[adressWuerzburg]{Christian Klingenberg}
\author[adressAs]{Jonas Kusch}

\address[adressWuerzburg]{University of Wuerzburg, Department of Mathematics,  Wuerzburg, Germany, \href{mailto:lena.baumann@uni-wuerzburg.de}{lena.baumann@uni-wuerzburg.de} (Lena Baumann), \href{mailto:klingen@mathematik.uni-wuerzburg.de}{christian.klingenberg@uni-wuerzburg.de} (Christian Klingenberg) }
\address[adressInnsbruck]{University of Innsbruck, Numerical Analysis and Scientific Computing, Innsbruck, Austria, \href{mailto:lukas.einkemmer@uibk.ac.at}{lukas.einkemmer@uibk.ac.at}}
\address[adressAs]{Norwegian University of Life Sciences, Department of Data Science, \r{A}s, Norway, \href{mailto:jonas.kusch@nmbu.no}{jonas.kusch@nmbu.no}}

\begin{abstract}
Computing numerical solutions of the thermal radiative transfer equations on a finely resolved grid can be costly due to high computational and memory requirements. A numerical reduced order method that has recently been applied to a wide variety of kinetic partial differential equations is the concept of dynamical low-rank approximation (DLRA). In this paper, we consider the thermal radiative transfer equations with Su-Olson closure, leading to a linearized kinetic model. For the conducted theoretical and practical considerations we use a multiplicative splitting of the distribution function that poses additional challenges in finding an energy stable discretization and deriving a hyperbolic Courant-Friedrichs-Lewy (CFL) condition. We propose such an energy stable DLRA scheme that makes use of the augmented basis update \& Galerkin integrator. This integrator allows for additional basis augmentations, enabling us to give a mathematically rigorous proof of energy stability and local mass conservation. Numerical examples confirm the derived properties and show the computational advantages of the DLRA scheme compared to a numerical solution of the full system of equations. 
\end{abstract}


\end{frontmatter}

\section{Introduction}\label{sec1:Introduction}

Thermal radiative transfer problems are a class of kinetic transport equations modeling the movement of particles that interact with a background medium, for instance by scattering or absorption. By this interaction, the background medium can heat up and itself emit new particles, enforcing the exchange of energy between particles and the background material. This process is described by two coupled equations, one for the particle density $f(t,x,\mu)$ and one for the internal energy $e(t,x)$ of the material. The variable $t \in \mathbb{R}^+$ denotes time, $x \in D \subset \mathbb{R}$ stands for the spatial and $\mu \in [-1,1]$ for the directional variable. The numerical solution of this set of equations is challenging as its high dimensionality requires enormous computational and memory costs. To overcome these problems, the concept of dynamical low-rank approximation (DLRA) \cite{koch2007dynamical} can be used. It is a numerical reduced order method providing accurate and efficient approximations of the solution of kinetic partial differential equations and has already been applied in various fields of research. Recent work for instance has been published on radiation transport \cite{Patwardhan2024,peng2020-2D,HauckSchnake2024}, radiation therapy \cite{stammerkusch2023}, plasma physics \cite{einkemmerlubich2019,einkemmerscalone2023}, chemical kinetics \cite{prugger2023} and Boltzmann type transport problems \cite{einkemmerhuying2021,huwang2022}. The main idea of DLRA consists in approximating the distribution function as
\begin{align*}
f(t,x,\mu) \approx \sum_{i,j=1}^r X_i(t,x) S_{ij}(t) V_j(t,\mu), 
\end{align*}
where $\{X_i : i=1,..,r \}$ are the orthonormal basis functions in space and $\{V_j : j=1,..,r \}$ the orthonormal basis functions in direction. This splitting can be understood as a continuous analogue to the singular value decomposition of a matrix, explaining why $\mathbf{S} = (S_{ij}) \in \mathbb{R}^{r \times r}$ is called the coefficient or coupling matrix and $r$ the rank of this approximation. However, note that the matrix $\mathbf{S}$ is not required to be diagonal. The time evolution of the low-rank factors is then determined by a projection of the equation onto the tangent space of the low-rank solution manifold. Integrators that ensure that the solution stays on the low-rank manifold while being robust to small singular values (otherwise this may lead to enormous restrictions on the choice of the time step size \cite{kieri2016discretized}) are the \textit{projector-splitting} \cite{lubich2014projector}, the \textit{(augmented) basis update \& Galerkin} (BUG) \cite{ceruti2022unconventional, ceruti2022rank}, and the \textit{parallel integrator} \cite{ceruti2024parallel}. 

A challenge for constructing a stable dynamical low-rank scheme is to derive a suitable discretization of the system. To account for the angular dependence of the solution, we consider a modal representation that makes use of a finite expansion of the particle density in terms of spherical harmonics ($P_N$). This approach is explained in \cite{bellglasstone1970,CaseZweifel1967} and for instance used in \cite{mcclarren2008-2,mcclarrenhauck2010,dargaville2019}. Also the spatial discretization has to be chosen carefully. In \cite{einkemmerhuying2021} it is shown that for an efficient dynamical low-rank scheme for the isothermal Boltzmann-BGK equation it is useful to consider a multiplicative splitting of the distribution function. To transfer knowledge to such systems, we decide on a multiplicative splitting of the particle density $f$ also for the considered problem. Different from \cite{baumann2024}, in this case it is per se not clear how to deal with the spatial derivatives. In addition, the multiplicative splitting poses further challenges for the proof of energy stability and the construction of a low-rank scheme to which we account for instance by pursuing a `first discretize, then low-rank' strategy. Recent work on this topic has been published for the linear Boltzmann-BGK equation in \cite{Baumann2024-BGK}. For the time discretization the potentially stiff opacity term has to be taken into account, leading to a coupled-implicit scheme similar to the one treated in \cite{baumann2024}, which again is complicated to solve.

In this paper we propose an energy stable multiplicative dynamical low-rank discretization for a linearized internal energy model called the Su-Olson problem. The main novelties are:
\begin{itemize}
    \item \textit{A multiplicative splitting of the distribution function:} Based on the insights from \cite{einkemmerhuying2021}, we consider a multiplicative splitting of the distribution function. Similar to \cite{Baumann2024-BGK}, we show that the spatial discretization has to be chosen carefully. Further, the multiplicative splitting requires additional modifications in the low-rank scheme in order to obtain an energy stable numerical discretization of the problem. 
    \item \textit{An energy stable numerical scheme with rigorous mathematical proofs:} We show that the derived DLRA scheme is energy stable and give rigorous mathematical proofs, enabling us to deduce a classic hyperbolic CFL condition. This allows to compute up to a maximal time step size of $\Delta t = \text{CFL} \cdot \Delta x$, which again enhances the performance of the algorithm. 
    \item \textit{A mass conservative augmented integrator:} We make use of the augmented BUG integrator from \cite{ceruti2022rank}, leading to a basis augmentation in two substeps of the numerical scheme. As this integrator allows for further modifications, we include additional basis augmentation steps that ensure the exactness of the projection operators needed for the theoretical proof of energy stability as well as the local conservation of mass.  
    \item \textit{Numerical test examples confirming the derived properties:} We compare the numerical results obtained from the low-rank discretization with the solution of the full system for relevant examples from the literature. The derived properties are validated and the accuracy and the efficiency of the proposed DLRA method is shown.
\end{itemize}
The structure of the paper is as follows: After the introduction in Section \ref{sec1:Introduction}, we provide background information on the thermal radiative transfer equations in Section \ref{sec2:Thermal radiative transfer}. We explain the considered multiplicative structure and derive two possible systems that in the continuous setting are equivalent. In Section \ref{sec3:Discretization}, we discretize both systems in angle, space and time. Section \ref{sec4:Energy stability} is devoted to the subject of energy stability. We show that the advection form of the multiplicative Su-Olson problem is generally not numerically stable in the sense of von Neumann, whereas for the conservative form a hyperbolic CFL condition is derived, under which an energy estimate can be given. Section \ref{sec5:DLRA} first provides an overview of the method of dynamical low-rank approximation. Then, a provably energy stable DLRA scheme is derived. The provided rigorous stability analysis is the main novelty of this work. In addition, mass conservation can be shown, when using a suitable truncation strategy. The numerical results in Section \ref{sec6:Numerical results} confirm the derived properties, before in Section \ref{sec7:Conclusion} a brief conclusion and outlook is given.

\section{Thermal radiative transfer}\label{sec2:Thermal radiative transfer}

In a one-dimensional setting the thermal radiative transfer equation with absorbing background material is given by
\begin{align*}
\frac{1}{c} \partial_t f(t,x,\mu) + \mu\partial_x f(t,x,\mu) &= \sigma(B(t,x)-f(t,x,\mu)),\\
\partial_t e(t,x) &= \sigma(\langle f(t,x,\cdot)-B(t,x)\rangle_{\mu}),
\end{align*}
where an integration over the directional domain $[-1,1]$ is denoted by $\langle \cdot \rangle_\mu$. The speed of light is denoted by $c$ and the variable $\sigma$ represents the opacity that encodes the rate at which particles are absorbed by the background medium. The black body radiation $B(T)$ at the material temperature $T$ can be described by the Stefan-Boltzmann law
\begin{align*}
    B(T) = acT^4,
\end{align*}
where $a = \frac{4\sigma_{\text{SB}}}{c}$ is the radiation density constant and $\sigma_{\text{SB}}$ the Stefan-Boltzmann constant. Further information on the thermal radiative transfer equations and their relevance in physics can be found in \cite{Pomraning1973,bellglasstone1970,HowellMengucDaunSiegel2020}. The above set of equations is not closed. To determine a relation between the temperature $T$ and the internal energy $e(T)$ we follow the ideas of Pomraning \cite{pomraning1979} and Su and Olson \cite{suolson1997} and set $e(T) = \alpha B(T)$. From this point on, we call $\alpha B(T)$ the internal energy of the material. Further, we perform a rescaling $\tau = \frac{t}{c}$ and by an abuse of notation write $t$ instead of $\tau$ in the remainder. This leads to the system
\begin{subequations}\label{eqs:Su-Olson}
\begin{align}
\partial_t f(t,x,\mu) + \mu\partial_x f(t,x,\mu) &= \sigma(B(t,x)-f(t,x,\mu))\label{eqs:Su-Olson f},\\
\partial_t B(t,x) &= \sigma(\langle f(t,x,\cdot)-B(t,x)\rangle_{\mu})\label{eqs:Su-Olson B},
\end{align}
\end{subequations}
where without loss of generality we assume $\alpha=1$. This system is a closed linearized internal energy model that is analytically solvable and commonly used as a benchmark for numerical examples \cite{olson2000,mcclarren2008-2, mcclarrenhauck2010}. In the following, we call equations \eqref{eqs:Su-Olson} the \textit{Su-Olson problem}. Note that for the moment we omit initial and boundary conditions.

In \cite{einkemmerhuying2021} it has been shown that for deriving an efficient dynamical low-rank scheme for the isothermal Boltzmann-BGK equation it is crucial to consider a multiplicative splitting of the distribution function. This allows one to separate a generally not low-rank Maxwellian from a remaining low-rank function $g$, to which the DLRA scheme is subsequently applied. The considered Su-Olson problem is similar in structure to this equation. To transfer knowledge of the construction of stable efficient dynamical low-rank schemes from the Su-Olson problem to more general kinetic equations, we have decided on a multiplicative splitting of the distribution function of the form
\begin{align}\label{eq:multiplicative splitting f=Bg}
f(t,x,\mu) = B(t,x) g(t,x,\mu),
\end{align}
and apply the low-rank ansatz to $g$. For this system, we give a mathematically rigorous proof of energy stability and derive a hyperbolic CFL condition, which allows for determining an optimal time step size, making the proposed scheme even more efficient. In this sense, this paper can be understood as an intermediate step from the Su-Olson problem treated in \cite{baumann2024} towards more complicated non-linear problems with multiplicative splitting. The main novelty of this paper consists in providing a rigorous theoretical analysis for the multiplicative Su-Olson problem, distinguishing this work for instance from the multiplicative DLRA scheme proposed in \cite{einkemmerhuying2021}, where the time step size of the constructed algorithm is chosen experimentally based on the corresponding numerical experiment. The idea of applying a multiplicative splitting to a linearized model and deriving a concrete CFL condition is also pursued in \cite{Baumann2024-BGK} for the linear Boltzmann-BGK equation. However, there are crucial differences between the multiplicative Su-Olson problem and the multiplicative linear Boltzmann-BGK equation. In the latter, due to the appearance of the Maxwellian equilibrium distribution, the evolution equation couples with the momentum equations. In the Su-Olson problem, the evolution equation for the particle density \eqref{eqs:Su-Olson f} couples through \eqref{eqs:Su-Olson B} with the background material. This leads to a different notion of stability in the theoretical analysis. In addition, the construction of the low-rank scheme requires new ideas for the basis augmentations as well as an adjusted truncation strategy to ensure mass conservation, which is not available in \cite{Baumann2024-BGK}. Hence, this paper provides novel insights into the construction of conservative and provably stable multiplicative low-rank schemes.

We insert ansatz \eqref{eq:multiplicative splitting f=Bg} into \eqref{eqs:Su-Olson f} and \eqref{eqs:Su-Olson B} and obtain the set of equations
\begin{subequations}
\label{eqs-naive: Bg full}
\begin{align}
\label{eq-naive: Bg full eq g}
\partial_t g (t,x,\mu) =& - \mu \partial_x g(t,x,\mu)  +\sigma \left(1-g(t,x,\mu) \right) - \frac{g(t,x,\mu)}{B(t,x)}\partial_t B(t,x)- \mu \frac{g(t,x,\mu)}{B(t,x)}\partial_x B(t,x)\\
\partial_t B(t,x) &= \sigma B(t,x)\left(\langle g(t,x,\mu) \rangle_{\mu}- 2 \right),\label{eq-naive: Bg full eq B}
\end{align}
\end{subequations}
that is called the \textit{advection form} of the multiplicative system. Using the product rule, it splits up the spatial derivatives for $B$ and $g$ in \eqref{eq-naive: Bg full eq g}. This corresponds to the form in which the multiplicative splitting in \cite{einkemmerhuying2021} is applied to the Boltzmann-BGK equation. Equation \eqref{eq-naive: Bg full eq g} can be equivalently rewritten into a \textit{conservative form,} leaving the spatial derivative of $Bg$ together and leading to the system
\begin{subequations}
\label{eqs: Bg full}
\begin{align}
\label{eq: Bg full eq g}
\partial_t g (t,x,\mu) =& - \frac{\mu}{B(t,x)}\partial_x \left(B(t,x) g(t,x,\mu) \right) +\sigma \left(1-g(t,x,\mu) \right)- \frac{g(t,x,\mu)}{B(t,x)}\partial_t B(t,x),\\
\label{eq: Bg full eq B}
\partial_t B(t,x) &= \sigma B(t,x)\left(\langle g(t,x,\mu) \rangle_{\mu}- 2 \right).
\end{align}
\end{subequations}

In later considerations, we are interested in the conservation properties of our numerical scheme. For the multiplicative Su-Olson problem, the mass and the momentum of the system shall be defined as follows.
\begin{definition}[Macroscopic quantities]\label{Def1: Macroscopic quantities}
The \textit{mass} of the multiplicative Su-Olson problem is defined as 
\begin{align*}
\rho (t,x) = \int f(t,x,\mu) \mathrm{d}\mu + B(t,x) = B(t,x) \int g(t,x,\mu) \mathrm{d}\mu + B(t,x).
\end{align*}
The \textit{momentum} is given as 
\begin{align*}
u(t,x) = \int \mu f(t,x,\mu) \mathrm{d}\mu = B(t,x) \int \mu g(t,x,\mu) \mathrm{d}\mu.
\end{align*}
\end{definition}
In particular, the multiplicative Su-Olson problem satisfies the local conservation law
\begin{align}\label{local conservation law}
\partial_t \rho (t,x) + \partial_x u(t,x) = 0. 
\end{align}
Global conservation of mass is then obtained by integrating over the spatial domain \cite{einkemmerlubich2019}. 

In the following, we discretize both sets of equations to compare them in terms of numerical stability. We derive an energy stable dynamical low-rank scheme and give a concrete hyperbolic CFL condition. Note that in contrast to \cite{baumann2024, einkemmerhuying2021}, but similar to \cite{Baumann2024-BGK}, we first discretize the equations and then apply the low-rank ansatz here.

\section{Discretization of the multiplicative system}\label{sec3:Discretization}

In this section, we fully discretize the advection form \eqref{eqs-naive: Bg full} as well as the conservative form \eqref{eqs: Bg full} of the multiplicative system. We start with the angular and spatial discretization, followed by the time discretization. 

\subsection{Angular discretization}

For the angular discretization a modal approach with normalized Legendre polynomials $P_\ell$ is used \cite{bellglasstone1970,CaseZweifel1967}. They constitute a complete set of orthogonal functions on the interval $[-1,1]$ that satisfy $\langle P_k, P_\ell \rangle_{\mu} = \delta_{k\ell}$. As an approximation, a finite expansion of the distribution function $g$ with $N_\mu$ expansion coefficients, called the moments, is used. It writes 
\begin{align*}
g(t,x,\mu)\approx g_{N_\mu}(t,x,\mu) = \sum_{\ell= 0}^{N_\mu-1} v_\ell(t,x) P_\ell(\mu).
\end{align*}
We insert this representation into \eqref{eqs-naive: Bg full}, multiply \eqref{eq-naive: Bg full eq g} with $P_k(\mu)$ and integrate over $\mu$.
Further, we introduce the matrix $\mathbf{A} \in \mathbb{R}^{N_\mu\times N_\mu}$ with entries $A_{k\ell} :=\langle P_k, \mu P_\ell \rangle_{\mu}$ and use that $P_0 = \frac{1}{\sqrt{2}}$. This gives
\begin{subequations}\label{eqs-naive: Bg angular}
\begin{align}
\label{eq-naive: Bg angular eq g} \partial_t v_k (t,x) =& - \sum_{\ell=0}^{N_\mu-1} \partial_x v_\ell(t,x) A_{k\ell} +\sigma \left(\sqrt{2} \delta_{k0} -v_k(t,x) \right)- \frac{v_k(t,x)}{B(t,x)} \partial_t B(t,x)\\
&- \sum_{\ell=1}^{N_\mu-1} \frac{v_\ell}{B(t,x)} \partial_x B(t,x) A_{k\ell},\nonumber\\
\partial_t B(t,x) =& \ \sigma B(t,x)\left(\sqrt{2} v_0(t,x) - 2 \right)\label{eq-naive: Bg angular eq B}.
\end{align}
\end{subequations}
Analogously, we obtain for system \eqref{eqs: Bg full} the angularly discretized equations
\begin{subequations}
\label{eqs: Bg angular}
\begin{align}
\label{eq: Bg angular eq g}
\partial_t v_k(t,x) =& - \sum_{\ell=0}^{N_\mu-1} \frac{1}{B(t,x)} \partial_x \left(B(t,x) v_\ell(t,x) \right) A_{k\ell} +\sigma \left( \sqrt{2} \delta_{k0}-v_k(t,x) \right)- \frac{v_k(t,x)}{B(t,x)}\partial_t B(t,x),\\
\label{eq: Bg angular eq B}
\partial_t B(t,x) =& \ \sigma B(t,x)\left(\sqrt{2} v_0(t,x) - 2 \right).
\end{align}
\end{subequations}
Note that the matrix $\mathbf{A}$ is symmetric and diagonalizable in the form $\mathbf{A} = \mathbf{Q}\mathbf{M}\mathbf{Q}^\top$ with $\mathbf{Q}$ orthonormal and $\mathbf{M} = \diag(\sigma_1,...,\sigma_{N_\mu})$. We then define $|\mathbf{A}| = \mathbf{Q} |\mathbf{M}| \mathbf{Q}^\top$. 

\subsection{Spatial discretization}

For the spatial discretization we prescribe a spatial grid with $N_x$ grid cells and equidistant spacing $\Delta x = \frac{1}{N_x}$. Spatially dependent quantities are approximated at the grid points $x_j$ for $j=1,...,N_x$, and denoted by
\begin{align*}
B_j(t) \approx B(t,x_j), \qquad v_{jk}(t) \approx v_k(t,x_j).
\end{align*}
First-order spatial derivatives $\partial_x$ are approximated using the tridiagonal stencil matrices $\mathbf{D}^x \in \mathbb{R}^{N_x \times N_x}$ and a second-order stabilization term $\mathbf{D}^{xx} \in \mathbb{R}^{N_x \times N_x}$ approximating $\frac{1}{2} \Delta x \partial_{xx}$ is added. The entries of those matrices are defined as 
\begin{align*}
D_{j,j\pm 1}^{x}= \frac{\pm 1}{2\Delta x}\;,\qquad D_{j,j}^{xx}= -\frac{1}{\Delta x}\;, \quad D_{j,j\pm 1}^{xx}= \frac{1}{2\Delta x}\;.
\end{align*}
Note that from now on we assume periodic boundary conditions to which we account by setting
\begin{align*}
D_{1,N_x}^{x} = \frac{-1}{2\Delta x}
\;,\quad D_{N_x,1}^{x} = \frac{1}{2\Delta x}\;, \quad D_{1,N_x}^{xx} = D_{N_x,1}^{xx} = \frac{1}{2 \Delta x}.
\end{align*}

The stencil matrices $\mathbf{D}^x$ and $\mathbf{D}^{xx}$ then fulfill the following properties:

\begin{lemma}\label{lemma:stencil matrices}
Let $y,z \in \mathbb{R}^{N_x}$ with indices $i,j=1,...,N_x$. 
It holds
\begin{align*}
\sum_{i,j=1}^{N_x} y_j D_{ji}^x z_i = - \sum_{i,j=1}^{N_x} z_j D_{ji}^{x} y_i\;, \hspace{0.15cm} \sum_{i,j=1}^{N_x} z_j D_{ji}^{x} z_i = 0 \;, \hspace{0.15cm} \sum_{i,j=1}^{N_x} y_j D_{ji}^{xx} z_i = \sum_{i,j=1}^{N_x} z_j D_{ji}^{xx} y_i.
\end{align*}
Moreover, let $\mathbf D^{+}\in\mathbb{R}^{N_x \times N_x}$ be defined as
\begin{align*}
D_{j,j}^{+}= \frac{- 1}{\sqrt{2\Delta x}}\;,\qquad D_{j,j + 1}^{+}= \frac{ 1}{\sqrt{2\Delta x}}\;.
\end{align*}
Then, $\sum_{i,j =1}^{N_x} z_j D_{ji}^{xx} z_i = -  \sum_{j=1}^{N_x} \left(\sum_{i=1}^{N_x} D_{ji}^+ z_i\right)^2$.
\end{lemma}
\begin{proof}
See \cite[Lemma 4.2]{baumann2024}.
\end{proof}

We insert the proposed discretization into the angularly discretized advection form \eqref{eqs-naive: Bg angular} of the equations and add a second-order stabilization term for $B \partial_x v$, leading to the angularly and spatially discretized set of equations
\begin{subequations}
\label{eqs-naive: Bg semidiscrete}
\begin{align}\label{eq-naive: Bg semidiscrete g}
\dot v_{jk}(t) =& - \sum_{i=1}^{N_x} \sum_{\ell=0}^{N_\mu-1} D_{ji}^x v_{i\ell}(t)  A_{k\ell} + \sum_{i=1}^{N_x} \sum_{\ell=0}^{N_\mu-1} D_{ji}^{xx} v_{i\ell}(t) |A|_{k\ell}\\\nonumber
&+\sigma \left(\sqrt{2} \delta_{k0} -v_{jk}(t) \right) - \frac{v_{jk}(t)}{B_j(t)} \dot B_j(t) - \sum_{i=1}^{N_x} \sum_{\ell=0}^{N_\mu-1} \frac{v_{j\ell}(t)}{B_j(t)} D_{ji}^x B_i(t) A_{k\ell},\\
\dot B_j(t) =& \ \sigma B_j(t)\left(\sqrt{2} v_{j0}(t) - 2 \right).\label{eq-naive: Bg semidiscrete B}
\end{align}
\end{subequations}
Inserting the discretization into the angularly discretized conservative form \eqref{eqs: Bg angular} of the equations and adding a second-order stabilization term to $\partial_x \left(Bv \right)$ gives
\begin{subequations}
\label{eqs: Bg semidiscrete}
\begin{align}
\label{eq: Bg semidiscrete g} \dot v_{jk}(t) =& - \sum_{i=1}^{N_x} \sum_{\ell=0}^{N_\mu-1} \frac{1}{B_j(t)} D_{ji}^x B_i(t) v_{i\ell}(t) A_{k\ell}+ \sum_{i=1}^{N_x} \sum_{\ell=0}^{N_\mu-1} \frac{1}{B_j(t)} D_{ji}^{xx} B_i(t) v_{i\ell}(t) |A|_{k\ell}\\
&+\sigma \left(\sqrt{2} \delta_{k0} -v_{jk}(t) \right) - \frac{v_{jk}(t)}{B_j(t)} \dot B_j(t),\nonumber\\
\dot B_j(t) =& \ \sigma B_j(t)\left(\sqrt{2} v_{j0}(t) - 2 \right).\label{eq: Bg semidiscrete B}
\end{align}
\end{subequations}

Note that due to the different structure of the equations the stabilization term in \eqref{eq-naive: Bg semidiscrete g} is applied to $B \partial_x v$, whereas in \eqref{eq: Bg semidiscrete g} it is added for $\partial_x \left(Bv \right)$.

\subsection{Time discretization}

From \cite{baumann2024} we know that constructing an energy stable scheme for the Su-Olson problem is challenging. For the advection form of the equations we start from \eqref{eqs-naive: Bg semidiscrete} and apply an explicit Euler step to transport terms. The potentially stiff absorption term is treated implicitly and the time derivative $\partial_t B$ is approximated by its difference quotient. We obtain the following fully discrete space-time discretization
\begin{subequations}
\label{eqs-naive: Bg fully discrete}
\begin{align}\label{eq-naive: Bg fully discrete g}
v_{jk}^1 =& \ v_{jk}^0 - \Delta t \sum_{i=1}^{N_x} \sum_{\ell=0}^{N_\mu-1} D_{ji}^x v_{i\ell}^0  A_{k\ell} + \Delta t \sum_{i=1}^{N_x} \sum_{\ell=0}^{N_\mu-1} D_{ji}^{xx} v_{i\ell}^0 |A|_{k\ell}\\\nonumber
&+\sigma \Delta t \left(\sqrt{2} \delta_{k0} -v_{jk}^1 \right) - \Delta t \frac{1}{B_j^0} \frac{B_j^1-B_j^0}{\Delta t} v_{jk}^1 - \Delta t \sum_{i=1}^{N_x} \sum_{\ell=0}^{N_\mu-1} \frac{v_{j\ell}^0}{B_j^0} D_{ji}^x B_i^0  A_{k\ell},\\
B_j^1 =& \ B_j^0 + \sigma \Delta t B_j^1 \left(\sqrt{2} v_{j0}^1 - 2 \right),\label{eq-naive: Bg fully discrete B}
\end{align}
\end{subequations}
that describes one time step from time $t_0$ to time $t_1 = t_0 + \Delta t$ and holds for all further time steps equivalently. For the conservative form \eqref{eqs: Bg semidiscrete} we again apply an explicit Euler step to the transport parts, treat the absorption terms implicitly and approximate $\partial_t B$ by its difference quotient. In addition, we add a factor $\frac{B^1}{B^0}$ in the absorption term of \eqref{eq: Bg semidiscrete g}. This gives the fully discrete scheme
\begin{subequations}
\label{eqs: Bg fully discrete}
\begin{align}\label{eq: Bg fully discrete g}
v_{jk}^1 =& \ v_{jk}^0 - \Delta t \sum_{i=1}^{N_x} \sum_{\ell=0}^{N_\mu-1} \frac{1}{B_j^0} D_{ji}^x B_i^0 v_{i\ell}^0  A_{k\ell} + \Delta t \sum_{i=1}^{N_x} \sum_{\ell=0}^{N_\mu-1} \frac{1}{B_j^0} D_{ji}^{xx} B_i^0 v_{i\ell}^0  |A|_{k\ell}\\\nonumber
&+\sigma \Delta t \frac{B_j^1}{B_j^0} \left(\sqrt{2} \delta_{k0} -v_{jk}^1 \right) - \Delta t \frac{1}{B_j^0} \frac{B_j^1-B_j^0}{\Delta t} v_{jk}^1,\\
B_j^1 =& \ B_j^0 + \sigma \Delta t B_j^1 \left(\sqrt{2} v_{j0}^1 - 2 \right).\label{eq: Bg fully discrete B}
\end{align}
\end{subequations}
Note that the evolution equations \eqref{eq-naive: Bg fully discrete B} and \eqref{eq: Bg fully discrete B} for the internal energy $B$ are the same in both schemes. The main difference of \eqref{eq-naive: Bg fully discrete g} and \eqref{eq: Bg fully discrete g} consists in the distinct second-order stabilization terms and the additional factor $\frac{B^1_j}{B^0_j}$ in \eqref{eq: Bg fully discrete g} that we will explain later when showing energy stability.

\section{Energy stability}\label{sec4:Energy stability}

The goal of this section is to investigate energy stability of the derived schemes. Note that this section is closely related to the considerations in \cite{baumann2024}. We first introduce the following notations.

\begin{definition}
In the following we write $u_{jk}^0 := B_j^0 v_{jk}^0$ and $u_{jk}^1 := B_j^1 v_{jk}^1$ at time $t_0$ and $t_1$, respectively. Note that $\mathbf{u}^0 = (u_{jk}^0) \in \mathbb{R}^{N_x \times N_\mu}$ corresponds to $f(t=0,x,\mu)$ and $\mathbf{v}^0 = (v_{jk}^0) \in \mathbb{R}^{N_x \times N_\mu}$ corresponds to $g(t=0,x,\mu)$ in \eqref{eq:multiplicative splitting f=Bg}.
\end{definition}

With this notation we can give the definition of the total energy of a fully discretized system. 

\begin{definition}[Total energy]
Let $\mathbf{u}^0 \in \mathbb{R}^{N_x \times N_\mu}$ be the fully discretized angular solution of the full Su-Olson problem and $\mathbf{B}^0 = (B_j^0) \in \mathbb{R}^{N_x}$ the internal energy at time $t_0$. Then, the \textit{total energy} at this time is defined as
\begin{align*}
E^0 := \frac{1}{2} \Vert \mathbf{u}^0 \Vert_F^2 + \frac{1}{2} \Vert \mathbf{B}^0 \Vert_E^2,  
\end{align*}
where $\Vert \cdot \Vert_F$ denotes the Frobenius and $\Vert \cdot \Vert_E$ the Euclidean norm. For $t_1 = t_0 + \Delta t$ this definition shall hold analogously. 
\end{definition}

\subsection{Advection form}

We start with the advection form \eqref{eqs-naive: Bg fully discrete} of the Su-Olson problem which is comparable to the considered low-rank discretization in \cite{einkemmerhuying2021} for the isothermal Boltzmann-BGK equation in the sense that the term $\partial_x \left(Bv \right)$ is split up into the sum of $B \partial_x v$ and $v \partial_x B$. We can show that this scheme is not, in general, von Neumann stable. 

\begin{theorem}\label{th:1}
There exist initial values $\mathbf{v}^0 \in \mathbb{R}^{N_x \times N_\mu}$ and $\mathbf{B}^0 \in \mathbb{R}^{N_x}$ such that the advection form \eqref{eqs-naive: Bg fully discrete} of the Su-Olson problem for $\sigma=0$ is not von Neumann stable.
\end{theorem}
\begin{proof}
Let us assume a solution $v_{jk}^0$ that is constant in space and direction, e.g. $v_{jk}^0 = 1$. For this solution all spatial derivatives are zero, i.e. the terms containing $\mathbf{D}^x \mathbf{v}^0$ and $\mathbf{D}^{xx} \mathbf{v}^0$ in \eqref{eq-naive: Bg fully discrete g} drop out. We further assume that for the opacity it holds $\sigma = 0$, i.e. the Su-Olson problem reduces to a simple advection equation. From \eqref{eq-naive: Bg fully discrete B} we thus obtain that $B_j^1 = B_j^0 = B_j$, i.e. the internal energy is constant in time. We insert these results into \eqref{eq-naive: Bg fully discrete g} and get
\begin{align*}
v_{jk}^1 = 1 - \Delta t \sum_{i=1}^{N_x} \sum_{\ell=0}^{N_\mu-1} \frac{1}{B_j} D_{ji}^x B_i  A_{k\ell}.
\end{align*}
Multiplication with $B_j$ then leads to
\begin{align*}
u_{jk}^1 = u_{jk}^0 - \Delta t \sum_{i=1}^{N_x} \sum_{\ell=0}^{N_\mu-1} D_{ji}^x u_{i\ell}^0  A_{k\ell}.
\end{align*}
This is a discretization of $\partial_t u + \mu \partial_x u = 0$ with an explicit Euler step forward in time and a centered finite difference scheme in space. From \cite{LeVeque2007} we know that this discretization is not von Neumann stable.
\end{proof}
This result matches our considerations that the DLRA scheme presented in \cite{einkemmerhuying2021} is not von Neumann stable but stable for relatively small time step sizes. Hence, Theorem \ref{th:1} is intended to serve as a motivation to seek a generally stable numerical discretization as done in the next section.

\subsection{Conservative form}

For the conservative form of the discretization of the Su-Olson problem given in \eqref{eqs: Bg fully discrete}, we can derive a hyperbolic CFL condition and show that under this time step restriction the total energy of the system decreases over time. We start with the following lemma.
\begin{lemma}\label{lemma:FourierAnalysis}
Under the time step restriction $\Delta t \leq \Delta x$ it holds
\begin{align*}
&\frac{\Delta t}{2} \sum_{i=1}^{N_x} \sum_{\ell=0}^{N_\mu-1} \left( \sum_{j=1}^{N_x} \sum_{k=0}^{N_\mu-1}  \left(D_{ji}^x u_{jk}^1 A_{k\ell} - D_{ji}^{xx} u_{jk}^1 |A|_{k\ell}\right)\right)^2 - \sum_{i=1}^{N_x} \sum_{\ell=0}^{N_\mu-1} \left( \sum_{j=1}^{N_x} \sum_{k=0}^{N_\mu-1} D_{ji}^+ u_{jk}^1 |A|_{k\ell}^{1/2} \right)^2 \leq 0.
\end{align*}
\end{lemma}
\begin{proof}
See \cite[Lemma 5.2]{baumann2024}.
\end{proof}
With this relation we can now show that the energy of the conservative form \eqref{eqs: Bg fully discrete} of the Su-Olson system dissipates and hence the system is energy stable.
\begin{theorem}\label{Theorem: Energy stability full system}
Under the time step restriction $\Delta t \leq \Delta x$ the fully discrete system \eqref{eqs: Bg fully discrete} is energy stable, i.e. it holds $E^1 \leq E^0$.
\end{theorem}
\begin{proof}
The proof of this theorem is similar to the proof of \cite[Theorem 5.3]{baumann2024}. We start with equation \eqref{eq: Bg fully discrete B} and multiply it with $B_j^1$. This gives
\begin{align*}
\left(B_j^1\right)^2 =& \ B_j^0 B_j^1 + \sigma \Delta t \left(B_j^1\right)^2 \left(\sqrt{2} v_{j0}^1 - 2 \right).
\end{align*}
Note that it holds
\begin{align*}
B_j^0 B_j^1 = \frac{1}{2} \left(B_j^1 \right)^2 + \frac{1}{2} \left(B_j^0 \right)^2 - \frac{1}{2} \left( B_j^1 - B_j^0\right)^2.
\end{align*}
We insert this relation and sum over $j$, giving
\begin{align}\label{eq: th1 proof B}
\frac{1}{2} \sum_{j=1}^{N_x} \left(B_j^1\right)^2 =& \ \frac{1}{2} \sum_{j=1}^{N_x} \left(B_j^0 \right)^2 - \frac{1}{2} \sum_{j=1}^{N_x} \left( B_j^1 - B_j^0\right)^2 + \sigma \Delta t \sum_{j=1}^{N_x} \left(B_j^1\right)^2 \left(\sqrt{2} v_{j0}^1 - 2 \right).
\end{align}
Next, we multiply equation \eqref{eq: Bg fully discrete g} with $B_j^1 B_j^0 v_{jk}^1$ and sum over $j$ and $k$. This leads to 
\begin{align}\label{eq: proof energy stability}
\sum_{j=1}^{N_x} \sum_{k=0}^{N_\mu-1} B_j^1 B_j^0 \left(v_{jk}^1\right)^2 =& \ \sum_{j=1}^{N_x} \sum_{k=0}^{N_\mu-1} B_j^0 v_{jk}^0 B_j^1 v_{jk}^1 - \Delta t \sum_{i,j=1}^{N_x} \sum_{k,\ell=0}^{N_\mu-1} B_j^1 v_{jk}^1 D_{ji}^x B_i^0 v_{i\ell}^0 A_{k\ell}\nonumber\\
&+ \Delta t \sum_{i,j=1}^{N_x} \sum_{k,\ell=0}^{N_\mu-1} B_j^1 v_{jk}^1 D_{ji}^{xx} B_i^0 v_{i\ell}^0 |A|_{k\ell} +\sigma \Delta t \sum_{j=1}^{N_x} \sum_{k=0}^{N_\mu-1} \left(B_j^1\right)^2 v_{jk}^1 \left(\sqrt{2} \delta_{k0} -v_{jk}^1 \right)\\
&- \sum_{j=1}^{N_x} \sum_{k=0}^{N_\mu-1} B_j^1 \left(v_{jk}^1\right)^2 \left(B_j^1-B_j^0 \right)\nonumber
\end{align}
Note that for this step the additional factor $\frac{B_j^1}{B_j^0}$ in the absorption term of \eqref{eq: Bg fully discrete g} is crucial. As above, it holds for the first term on the right-hand side that
\begin{align*}
&\sum_{j=1}^{N_x} \sum_{k=0}^{N_\mu-1} B_j^0 v_{jk}^0 B_j^1 v_{jk}^1
= \sum_{j=1}^{N_x} \sum_{k=0}^{N_\mu-1} \left( \frac{1}{2} \left( B_j^1 v_{jk}^1 \right)^2 + \frac{1}{2} \left( B_j^0 v_{jk}^0 \right)^2 - \frac{1}{2} \left( B_j^1 v_{jk}^1 - B_j^0 v_{jk}^0 \right)^2 \right).
\end{align*}
We insert the notation $u_{jk}^0 = B_j^0 v_{jk}^0$ and $u_{jk}^1 = B_j^1 v_{jk}^1$, respectively, as well as insert the above relation into \eqref{eq: proof energy stability}, bring the last term of \eqref{eq: proof energy stability} to the left-hand side and rearrange the equation. We obtain
\begin{align}\label{Theorem2-equ}
\frac{1}{2} \sum_{j=1}^{N_x} \sum_{k=0}^{N_\mu-1} \left(u_{jk}^1\right)^2 =& \ \frac{1}{2} \sum_{j=1}^{N_x} \sum_{k=0}^{N_\mu-1} \left( u_{jk}^0\right)^2 - \frac{1}{2} \sum_{j=1}^{N_x} \sum_{k=0}^{N_\mu-1} \left( u_{jk}^1 - u_{jk}^0\right)^2 - \Delta t \sum_{i,j=1}^{N_x} \sum_{k,\ell=0}^{N_\mu-1} u_{jk}^1 D_{ji}^x u_{i\ell}^0 A_{k\ell}\\
&+ \Delta t \sum_{i,j=1}^{N_x} \sum_{k,\ell=0}^{N_\mu-1}  u_{jk}^1 D_{ji}^{xx} u_{i\ell}^0 |A|_{k\ell} +\sigma \Delta t \sum_{j=1}^{N_x} \sum_{k=0}^{N_\mu-1} \left(B_j^1\right)^2 v_{jk}^1 \left(\sqrt{2} \delta_{k0} -v_{jk}^1 \right).\nonumber
\end{align}
We now add the zero term $\Delta t \sum_{i,j=1}^{N_x} \sum_{k,\ell=0}^{N_\mu-1} u_{jk}^1 D_{ji}^x u_{i\ell}^1 A_{k\ell}$ and add and subtract the second-order term $\Delta t \sum_{i,j=1}^{N_x} \sum_{k,\ell=0}^{N_\mu-1}  u_{jk}^1 D_{ji}^{xx} u_{i\ell}^1 |A|_{k\ell}$ giving
\begin{align*}
\frac{1}{2} \sum_{j=1}^{N_x} \sum_{k=0}^{N_\mu-1} \left(u_{jk}^1\right)^2 =& \ \frac{1}{2} \sum_{j=1}^{N_x} \sum_{k=0}^{N_\mu-1} \left( u_{jk}^0\right)^2 - \frac{1}{2} \sum_{j=1}^{N_x} \sum_{k=0}^{N_\mu-1} \left( u_{jk}^1 - u_{jk}^0\right)^2 - \Delta t \sum_{i,j=1}^{N_x} \sum_{k,\ell=0}^{N_\mu-1} u_{jk}^1 D_{ji}^x \left( u_{i\ell}^0 - u_{i\ell}^1 \right) A_{k\ell}\\
&+ \Delta t \sum_{i,j=1}^{N_x} \sum_{k,\ell=0}^{N_\mu-1}  u_{jk}^1 D_{ji}^{xx} \left( u_{i\ell}^0 - u_{i\ell}^1 \right) |A|_{k\ell} + \Delta t \sum_{i,j=1}^{N_x} \sum_{k,\ell=0}^{N_\mu-1}  u_{jk}^1 D_{ji}^{xx} u_{i\ell}^1 |A|_{k\ell}\nonumber\\
&+\sigma \Delta t \sum_{j=1}^{N_x} \sum_{k=0}^{N_\mu-1} \left(B_j^1\right)^2 v_{jk}^1 \left(\sqrt{2} \delta_{k0} -v_{jk}^1 \right).
\end{align*}
In the next step, we apply Young's inequality, which states that for $a,b\in\mathbb{R}$ we have $a\cdot b \leq \frac{a^2}{2} + \frac{b^2}{2}$, to the term
\begin{align*}
&- \Delta t \sum_{i,j=1}^{N_x} \sum_{k,\ell=0}^{N_\mu-1} u_{jk}^1 D_{ji}^x \left( u_{i\ell}^0 - u_{i\ell}^1 \right) A_{k\ell} + \Delta t \sum_{i,j=1}^{N_x} \sum_{k,\ell=0}^{N_\mu-1}  u_{jk}^1 D_{ji}^{xx} \left( u_{i\ell}^0 - u_{i\ell}^1 \right) |A|_{k\ell}\\
=& - \Delta t \sum_{i=1}^{N_x} \sum_{\ell=0}^{N_\mu-1} \left( u_{i\ell}^0 - u_{i\ell}^1 \right) \left( \sum_{j=1}^{N_x} \sum_{k=0}^{N_\mu-1}  \left(D_{ji}^x u_{jk}^1 A_{k\ell} - D_{ji}^{xx} u_{jk}^1 |A|_{k\ell}\right)\right)\\
\leq& \ \frac{1}{2} \sum_{i=1}^{N_x} \sum_{\ell=0}^{N_\mu-1} \left( u_{i\ell}^0 - u_{i\ell}^1 \right)^2 + \frac{\left( \Delta t \right)^2}{2} \sum_{i=1}^{N_x} \sum_{\ell=0}^{N_\mu-1} \left( \sum_{j=1}^{N_x} \sum_{k=0}^{N_\mu-1}  \left(D_{ji}^x u_{jk}^1 A_{k\ell} - D_{ji}^{xx} u_{jk}^1 |A|_{k\ell}\right)\right)^2.
\end{align*}
In addition, we note that with the properties of the stencil matrices from Lemma \ref{lemma:stencil matrices} we can write 
\begin{align*}
\Delta t \sum_{i,j=1}^{N_x} \sum_{k,\ell=0}^{N_\mu-1}  u_{jk}^1 D_{ji}^{xx} u_{i\ell}^1 |A|_{k\ell} = -\Delta t \sum_{i=1}^{N_x} \sum_{\ell=0}^{N_\mu-1} \left( \sum_{j=1}^{N_x} \sum_{k=0}^{N_\mu-1} D_{ji}^+ u_{jk}^1 |A|_{k\ell}^{1/2} \right)^2.
\end{align*}
We insert both relations and get
\begin{align*}
\frac{1}{2} \sum_{j=1}^{N_x} \sum_{k=0}^{N_\mu-1} \left(u_{jk}^1\right)^2 \leq& \ \frac{1}{2} \sum_{j=1}^{N_x} \sum_{k=0}^{N_\mu-1} \left( u_{jk}^0\right)^2 + \frac{\left( \Delta t \right)^2}{2} \sum_{i=1}^{N_x} \sum_{\ell=0}^{N_\mu-1} \left( \sum_{j=1}^{N_x} \sum_{k=0}^{N_\mu-1}  \left(D_{ji}^x u_{jk}^1 A_{k\ell} - D_{ji}^{xx} u_{jk}^1 |A|_{k\ell}\right)\right)^2\\
&- \Delta t \sum_{i=1}^{N_x} \sum_{\ell=0}^{N_\mu-1} \left( \sum_{j=1}^{N_x} \sum_{k=0}^{N_\mu-1} D_{ji}^+ u_{jk}^1 |A|_{k\ell}^{1/2} \right)^2 + \sigma \Delta t \sum_{j=1}^{N_x} \sum_{k=0}^{N_\mu-1} \left(B_j^1\right)^2 v_{jk}^1 \left(\sqrt{2} \delta_{k0} -v_{jk}^1 \right).
\end{align*}
With Lemma \ref{lemma:FourierAnalysis} we have that under the time step restriction $\Delta t \leq \Delta x$ it holds
\begin{align}\label{eq: th1 proof u}
\frac{1}{2} \sum_{j=1}^{N_x} \sum_{k=0}^{N_\mu-1} \left(u_{jk}^1\right)^2 \leq \frac{1}{2} \sum_{j=1}^{N_x} \sum_{k=0}^{N_\mu-1} \left( u_{jk}^0\right)^2 + \sigma \Delta t \sum_{j=1}^{N_x} \sum_{k=0}^{N_\mu-1} \left(B_j^1\right)^2 v_{jk}^1 \left(\sqrt{2} \delta_{k0} -v_{jk}^1 \right).
\end{align}
To obtain an expression for the total energy of the system, we add \eqref{eq: th1 proof u} and \eqref{eq: th1 proof B}. This gives
\begin{align*}
E^1 \leq& \ E^0  - \frac{1}{2} \sum_{j=1}^{N_x} \left( B_j^1 - B_j^0\right)^2 + \sigma \Delta t \sum_{j=1}^{N_x} \sum_{k=0}^{N_\mu-1} \left(B_j^1\right)^2 v_{jk}^1 \left(\sqrt{2} \delta_{k0} -v_{jk}^1 \right) + \sigma \Delta t \sum_{j=1}^{N_x} \left(B_j^1\right)^2 \left(\sqrt{2} v_{j0}^1 - 2 \right).
\end{align*}
The term $- \frac{1}{2} \sum_{j=1}^{N_x} \left( B_j^1 - B_j^0\right)^2$ is non-positive. The remaining two terms on the right-hand side can be rewritten and bounded as follows:
\begin{align*}
&\ \sigma \Delta t \sum_{j=1}^{N_x} \sum_{k=0}^{N_\mu-1} \left(B_j^1\right)^2 v_{jk}^1 \left(\sqrt{2} \delta_{k0} -v_{jk}^1 \right) + \sigma \Delta t \sum_{j=1}^{N_x} \left(B_j^1\right)^2 \left(\sqrt{2} v_{j0}^1 - 2 \right)\\
\leq& \ \sigma \Delta t \sum_{j=1}^{N_x} \sum_{k=0}^{N_\mu-1} \left(B_j^1\right)^2  \left(-\left(v_{jk}^1\right)^2 + 2\sqrt{2} v_{jk}^1 \delta_{k0} -2 \delta_{k0} \right)\\
=& - \sigma \Delta t \sum_{j=1}^{N_x} \sum_{k=0}^{N_\mu-1} \left(B_j^1\right)^2  \left( v_{jk}^1 - \sqrt{2} \delta_{k0} \right)^2 \leq 0.
\end{align*}
Hence, we have shown that under the time step restriction $\Delta t \leq \Delta x$ it holds $E^1 \leq E^0$, and the system is energy stable.
\end{proof}

\section{Dynamical low-rank approximation for the energy stable system}\label{sec5:DLRA}

Having attained an energy stable discretization of the multiplicative Su-Olson problem, its practical implementation can still pose numerical challenges such as large memory demands and computational costs, especially in a higher-dimensional setting. To overcome these problems, we introduce the concept of dynamical low-rank approximation.

\subsection{Background on DLRA}

The method of dynamical low-rank approximation has originally been introduced in a semi-discrete time-continuous matrix setting \cite{koch2007dynamical}, in which it will also be explained in this section. Let us consider the matrix differential equation
\begin{align*}
\dot{\mathbf{f}}(t) = \mathbf{F} \left(\mathbf{f}(t) \right),
\end{align*}
where $\mathbf{f}(t) \in \mathbb{R}^{N_x \times N_\mu}$ is the solution of the equation and $\mathbf{F} \left(\mathbf{f}(t) \right): \mathbb{R}^{N_x \times N_\mu} \to \mathbb{R}^{N_x \times N_\mu}$ denotes its right-hand side. The dynamical low-rank approximation of $\mathbf{f}(t)$ is then given by
\begin{align}\label{eq:low-rank f general form}
\mathbf{f}_r(t) = \mathbf{X}(t) \mathbf{S}(t) \mathbf{V}(t)^\top,
\end{align}
where $\mathbf{X}(t) \in \mathbf{R}^{N_x \times r}$ is the orthonormal basis in space, $\mathbf{V}(t) \in \mathbf{R}^{N_\mu \times r}$ the orthonormal basis in direction, and $\mathbf{S}(t) \in \mathbb{R}^{r \times r}$ the coupling or coefficient matrix of the rank $r$ approximation. All matrices of the form \eqref{eq:low-rank f general form} constitute the manifold of low-rank matrices $\mathcal{M}_r$. The basis matrices $\mathbf{X}(t)$ and $\mathbf{V}(t)$, and the coefficient matrix $\mathbf{S}(t)$ shall then be evolved in time such that the minimization problem
\begin{align*}
\min_{\dot{\mathbf{f}}_r(t) \in \mathcal{T}_{\mathbf{f}_r(t)} \mathcal{M}_r} \Vert \dot{\mathbf{f}}_r(t) - \mathbf{F}\left(\mathbf{f}_r(t)\right)\Vert_F 
\end{align*}
is fulfilled at all times $t$, where $\mathcal{T}_{\mathbf{f}_r(t)} \mathcal{M}_r$ denotes the tangent space of the low-rank manifold $\mathcal{M}_r$ at $\mathbf{f}_r(t)$. In \cite{koch2007dynamical}, it has been shown that solving this minimization problem is equivalent to projecting the right-hand side $\mathbf{F}\left(\mathbf{f}_r(t)\right)$ by means of an orthogonal projection $\mathbf{P}$ onto $\mathcal{T}_{\mathbf{f}_r(t)} \mathcal{M}_r$, and solving 
\begin{align}\label{sec5: DLRA projection}
\dot{\mathbf{f}}_r(t) = \mathbf{P} \left(\mathbf{f}_r(t) \right) \mathbf{F}\left(\mathbf{f}_r(t)\right).
\end{align}
For $\mathbf{f}_r = \mathbf{X}\mathbf{S} \mathbf{V}^\top$, the orthogonal projection $\mathbf{P}$ onto $\mathcal{T}_{\mathbf{f}_r(t)} \mathcal{M}_r$ is explicitly given \cite{koch2007dynamical} as 
\begin{align*}
\mathbf{P}(\mathbf{f}_r) \mathbf{F}\left(\mathbf{f}_r\right) = \mathbf{X} \mathbf{X}^\top \mathbf{F}\left(\mathbf{f}_r\right) - \mathbf{X} \mathbf{X}^\top \mathbf{F}\left(\mathbf{f}_r\right) \mathbf{V} \mathbf{V}^\top + \mathbf{F}\left(\mathbf{f}_r\right) \mathbf{V} \mathbf{V}^\top.
\end{align*}
Different time integrators that make use of the special form of this orthogonal projection such as the projector-splitting \cite{lubich2014projector}, the (augmented) basis update \& Galerkin (BUG) \cite{ceruti2022unconventional, ceruti2022rank}, and the parallel integrator \cite{ceruti2024parallel} exist. They are able to evolve the solution on the low-rank manifold while not suffering from the stiffness of \eqref{sec5: DLRA projection}. In this paper, we focus on the augmented basis update \& Galerkin (BUG) integrator, that shall be explained in the following.

The BUG integrator first updates and augments the spatial basis $\mathbf{X}$ and the directional basis $\mathbf{V}$ in parallel. This leads to an increase from rank $r$ to rank $2r$. Note that we denote augmented quantities of rank $2r$ with hats. Next, a Galerkin step is conducted for the coefficient matrix $\mathbf{S}$ in the augmented setting, before in a last step all augmented quantities are truncated to a new rank $r_1 \leq 2r$. To be more specific, the BUG integrator evolves the low-rank solution $\mathbf{f}_r^0 = \mathbf{X}^0 \mathbf{S}^0 \mathbf{V}^{0,\top}$ at time $t_0$ to the time-updated low-rank solution $\mathbf{f}_r^1 = \mathbf{X}^1 \mathbf{S}^1 \mathbf{V}^{1,\top}$ at time $t_1 = t_0 + \Delta t$ as follows:

\textbf{\textit{K}-Step}:
We fix the directional basis $\mathbf{V}^0$ at time $t_0$ and introduce the notation $\mathbf{K}(t) = \mathbf{X}(t) \mathbf{S}(t)$. Then, we update the spatial basis from $\mathbf{X}^0$ to $\widehat{\mathbf{X}}^1 \in \mathbb{R}^{N_x \times 2r}$ by solving the PDE
\begin{align*}
\dot{\mathbf{K}}(t) = \mathbf{F}\left(\mathbf{K}(t) \mathbf{V}^{0,\top} \right) \mathbf{V}^0, \quad \mathbf{K}(t_0) = \mathbf{X}^0 \mathbf{S}^0,
\end{align*}
and then determining $\widehat{\mathbf{X}}^1$ as an orthonormal basis of the augmented matrix $[\mathbf{K}(t_1), \mathbf{X}^0] \in \mathbb{R}^{N_x \times 2r}$, e.g. by QR-decomposition. We store $\widehat{\mathbf{M}} = \widehat{\mathbf{X}}^{1,\top} \mathbf{X}^0 \in \mathbb{R}^{2r \times r}$.

\textbf{\textit{L}-Step}:
We fix the spatial basis $\mathbf{X}^0$ at time $t_0$ and introduce the notation $\mathbf{L}(t) = \mathbf{V}(t) \mathbf{S}(t)^\top$. Then, we update the directional basis from $\mathbf{V}^0$ to $\widehat{\mathbf{V}}^1 \in \mathbb{R}^{N_\mu \times 2r}$ by solving the PDE
\begin{align*}
\dot{\mathbf{L}}(t) = \mathbf{F}\left( \mathbf{X}^0 \mathbf{L}(t)^\top \right)^\top \mathbf{X}^0, \quad \mathbf{L}(t_0) = \mathbf{V}^0 \mathbf{S}^{0,\top},
\end{align*}
and then determining $\widehat{\mathbf{V}}^1$ as an orthonormal basis of the augmented matrix $[\mathbf{L}(t_1), \mathbf{V}^0] \in \mathbb{R}^{N_\mu \times 2r}$, e.g. by QR-decomposition. We store $\widehat{\mathbf{N}} = \widehat{\mathbf{V}}^{1,\top} \mathbf{V}^0 \in \mathbb{R}^{2r \times r}$. 

\textbf{\textit{S}-step}: We update the coupling matrix from $\mathbf{S}^0 \in \mathbb{R}^{r \times r}$ to $\widehat{\mathbf{S}}^1 \in \mathbb{R}^{2r \times 2r}$ by solving the ODE
\begin{align*}
\dot{\widehat{\mathbf{S}}}(t) = \widehat{\mathbf{X}}^{1, \top} \mathbf{F} \left( \widehat{\mathbf{X}}^1 \widehat{\mathbf{S}}(t) \widehat{\mathbf{V}}^{1,\top} \right) \widehat{\mathbf{V}}^1, \quad \widehat{\mathbf{S}}(t_0) = \widehat{\mathbf{M}} \mathbf{S}^0 \widehat{\mathbf{N}}^\top.
\end{align*}

\textbf{Truncation}: We compute the singular value decomposition of $\widehat{\mathbf{S}}^1 = \widehat{\mathbf{P}} \mathbf{\Sigma}\widehat {\mathbf{Q}}^\top$ with $\mathbf{\Sigma} = \text{diag}(\sigma_j)$. The new rank $r_1 \leq 2r$ is chosen such that for a prescribed tolerance parameter $\vartheta$ it holds
\begin{align*}
\left(\sum_{j=r_1+1}^{2r} \sigma_j^2\right)^{1/2} \leq \vartheta.
\end{align*}
We set $\mathbf{S}^1 \in \mathbb{R}^{r_1\times r_1}$ to be the matrix containing the $r_1$ largest singular values. For the update of the spatial and the directional basis we introduce the matrices $\mathbf{P}^1 \in \mathbb{R}^{2r\times r_1}$ and $\mathbf{Q}^1 \in \mathbb{R}^{2r\times r_1}$ containing the first $r_1$ columns of $\widehat{\mathbf{P}}$ and $\widehat {\mathbf{Q}}$, respectively, and set $\mathbf{X}^1 = \widehat{\mathbf{X}}^1\mathbf{P}^1 \in \mathbb{R}^{N_x \times r_1}$ and $\mathbf{V}^1 = \widehat{\mathbf{V}}^1 \mathbf{Q}^1 \in \mathbb{R}^{N_\mu \times r_1}$.

Altogether, this gives the time-updated low-rank approximation in the form $\mathbf{f}_r^1 = \mathbf{X}^1 \mathbf{S}^1 \mathbf{V}^{1,\top}$ at time $t_1 = t_0 + \Delta t$. Note that in the following, in an abuse of notation, we will write $\mathbf{f}$ instead of $\mathbf{f}_r$ and call this the low-rank solution of the considered problem.

\subsection{DLRA scheme for multiplicative Su-Olson}

In this section, the DLRA method is applied to the energy stable conservative form \eqref{eqs: Bg fully discrete} of the Su-Olson problem to evolve $\mathbf{v}^0 = \left(v_{jk}^0\right)$ to $\mathbf{v}^1 = \left( v_{jk}^1 \right)$. First note that for the derivation of the low-rank scheme we rewrite equations \eqref{eqs: Bg fully discrete}. In \eqref{eq: Bg fully discrete g}, we put all terms containing $v_{jk}^1$ to the left-hand side and divide by $1+\sigma \Delta t$. Further, we multiply \eqref{eq: Bg fully discrete B} with $\frac{1}{B_j^0}$. This establishes the system
\begin{subequations}\label{eqs: full reformulation for DLRA both}
\begin{align}
\frac{B_j^1}{B_j^0} v_{jk}^1 =& \ \frac{1}{1 + \sigma \Delta t }
v_{jk}^0 - \frac{\Delta t}{1 + \sigma \Delta t} \sum_{i=1}^{N_x} \sum_{\ell=0}^{N_\mu-1} \frac{1}{B_j^0} D_{ji}^x B_i^0 v_{i\ell}^0  A_{k\ell}\label{eqs: full reformulation for DLRA f}\\
&+ \frac{\Delta t}{1 + \sigma \Delta t} \sum_{i=1}^{N_x} \sum_{\ell=0}^{N_\mu-1} \frac{1}{B_j^0} D_{ji}^{xx} B_i^0 v_{i\ell}^0  |A|_{k\ell} +\frac{\sqrt{2} \sigma \Delta t}{1 + \sigma \Delta t} \frac{B_j^1}{B_j^0} \delta_{k0},\nonumber\\
\frac{B_j^1}{B_j^0} =& \ 1 + \sigma \Delta t \frac{B_j^1}{B_j^0} \left(\sqrt{2} v_{j0}^1 - 2 \right).\label{eqs: full reformulation for DLRA B}
\end{align}
\end{subequations}
We then apply the DLRA approach to this set of equations. In a first step, we want to update $v_{jk}^0 = \sum_{m,n=1}^r X_{jm}^0 S_{mn}^0 V_{kn}^0$ to $\frac{B_j^1}{B_j^0} v_{jk}^\ast = \sum_{m,n=1}^{4r} \doublehat{X}_{jm}^\ast \doublehat{S}_{mn}^\ast \doublehat{V}_{kn}^\ast$ for $k \neq 0$. This corresponds to an update of the advection part of \eqref{eqs: full reformulation for DLRA f} without the absorption term coupling to equation \eqref{eqs: full reformulation for DLRA B}. We introduce the notation $K_{jn}^0 = \sum_{m=1}^r X_{jm}^0 S_{mn}^0$ and solve the \textit{K}-step equation
\begin{subequations}\label{DLRA-scheme}
\begin{align}\label{DLRA-scheme: K step}
K_{jp}^\ast =& \ \frac{1}{1 + \sigma \Delta t }
K_{jp}^0 - \frac{\Delta t}{1 + \sigma \Delta t} \frac{1}{B_j^0} \sum_{i=1}^{N_x} D_{ji}^x B_i^0 \sum_{n=1}^r K_{in}^0 \sum_{k,\ell=0}^{N_\mu-1} V_{\ell n}^0  A_{k\ell} V_{kp}^0\\
&+ \frac{\Delta t}{1 + \sigma \Delta t} \frac{1}{B_j^0} \sum_{i=1}^{N_x} D_{ji}^{xx} B_i^0 \sum_{n=1}^r K_{in}^0 \sum_{k,\ell=0}^{N_\mu-1}  V_{\ell n}^0  |A|_{k\ell} V_{kp}^0.\nonumber
\end{align}
We derive the updated basis $\widehat{\mathbf{X}}^\ast$ of rank $2r$ from $\widehat{\mathbf{X}}^\ast = \text{qr}\left([\mathbf{K}^\ast, \mathbf{X}^0]\right)$. Moreover, we perform an additional basis augmentation step according to
\begin{align}\label{DLRA-scheme: X basis augmentation 3r}
\doublehat{\mathbf{X}}^\ast = \text{qr} \left(\left[\widehat{\mathbf{X}}^\ast, \frac{1}{\mathbf{B}^0} \odot \mathbf{D}^x \left( \mathbf{B}^0 \odot \mathbf{X}^0 \right), \frac{1}{\mathbf{B}^0} \odot \mathbf{D}^{xx} \left( \mathbf{B}^0 \odot \mathbf{X}^0 \right)\right] \right).
\end{align}
This basis augmentation step is crucial for the theoretical proof of energy stability of the DLRA scheme as it ensures the exactness of the corresponding projection operators. Its concrete form is motivated from the theoretical stability analysis and will become clear in the proof of energy stability in the next section. Here, the symbol $\odot$ stands for a pointwise multiplication and the vector $\frac{1}{\mathbf{B}^0} \in \mathbb{R}^{N_x}$ is defined to contain the element $\frac{1}{B_j^0}$ for each $j=1,...,N_x$. In addition, we compute and store  $\doublehat{\mathbf{M}} = \doublehat{\mathbf{X}}^{\ast,\top} \mathbf{X}^0$. Note that in this paper we perform full rank updates, leading to an increase from rank $2r$ to $4r$. Quantities of rank $2r$ are denoted with one single hat and quantities of rank $4r$ with double hats. 

The \textit{L}-step can be computed in parallel with the \textit{K}-step. We introduce the notation $L_{mk}^0 = \sum_{n=1}^r S_{nm}^0 V_{nk}^0$ and solve
\begin{align}\label{DLRA-scheme: L step}
L_{kp}^\ast =& \ \frac{1}{1 + \sigma \Delta t}
L_{kp}^0 - \frac{\Delta t}{1 + \sigma \Delta t} \sum_{\ell=0}^{N_\mu-1} \sum_{m=1}^r A_{\ell k} L_{\ell m}^0 \sum_{i=1}^{N_x} X_{im}^0 B_i^0 \sum_{j=1}^{N_x} D_{ij}^x \frac{1}{B_j^0} X_{jp}^0\\
&+ \frac{\Delta t}{1 + \sigma \Delta t} \sum_{\ell=0}^{N_\mu-1} \sum_{m=1}^r |A|_{\ell k} L_{\ell m}^0 \sum_{i=1}^{N_x} X_{im}^0 B_i^0 \sum_{j=1}^{N_x} D_{ij}^{xx} \frac{1}{B_j^0} X_{jp}^0.\nonumber
\end{align}
We derive the updated basis $\widehat{\mathbf{V}}^\ast$ of rank $2r$ from $\widehat{\mathbf{V}}^\ast = \text{qr}\left([\mathbf{L}^\ast, \mathbf{V}^0]\right)$. 
Moreover, we perform an additional basis augmentation step according to
\begin{align}\label{DLRA-scheme: V basis augmentation 3r}
\doublehat{\mathbf{V}}^\ast = \text{qr} \left(\left[\widehat{\mathbf{V}}^\ast, \mathbf{A}^\top \mathbf{V}^0, \left|\mathbf{A}\right|^\top \mathbf{V}^0 \right] \right).
\end{align}
Again, this additional basis augmentation step is required for the exactness of the corresponding projection operators in the theoretical proof of energy stability. Its concrete form will become clear in the next section. In addition, we compute and store $\doublehat{\mathbf{N}} = \doublehat{\mathbf{V}}^{\ast,\top} \mathbf{V}^0$. 

For the \textit{S}-step we use the information from the \textit{K}- and \textit{L}-step, set $\widetilde{S}_{mn}^0 = \sum_{j,k=1}^r \doublehat{M}_{mj} S_{jk}^0 \doublehat{N}_{nk}$, and solve
\begin{align}\label{DLRA-scheme: S step}
\doublehat{S}_{qp}^\ast =& \ \frac{1}{1 + \sigma \Delta t} \widetilde{S}_{qp}^0 - \frac{\Delta t}{1 + \sigma \Delta t} \sum_{j=1}^{N_x} \doublehat{X}_{jq}^\ast \frac{1}{B_j^0} \sum_{i=1}^{N_x} D_{ji}^x B_i^0 \sum_{m,n=1}^{4r} \doublehat{X}_{im}^\ast \widetilde{S}_{mn}^0 \sum_{k,\ell=0}^{N_\mu-1} \doublehat{V}_{\ell n}^\ast A_{k\ell} \doublehat{V}_{kp}^\ast\\
&+ \frac{\Delta t}{1 + \sigma \Delta t} \sum_{j=1}^{N_x} \doublehat{X}_{jq}^\ast \frac{1}{B_j^0} \sum_{i=1}^{N_x} D_{ji}^{xx} B_i^0 \sum_{m,n=1}^{4r} \doublehat{X}_{im}^\ast \widetilde{S}_{mn}^0 \sum_{k,\ell=0}^{N_\mu-1} \doublehat{V}_{\ell n}^\ast |A|_{k\ell} \doublehat{V}_{kp}^\ast\nonumber.
\end{align}

In the next step, we consider the equations for $k=0$. In this case, the equations for $\frac{B_j^1}{B_j^0} v_{j0}^1$ and $\frac{B_j^1}{B_j^0}$ couple and we solve the system
\begin{align}
\frac{B_j^1}{B_j^0} v_{j0}^1 =& \ \frac{1}{1 + \sigma \Delta t} \sum_{m,n=1}^r X_{jm}^0 S_{mn}^0 V_{0n}^0 - \frac{\Delta t}{1 + \sigma \Delta t} \sum_{i=1}^{N_x} \frac{1}{B_j^0} D_{ji}^x B_i^0 \sum_{m,n=1}^{4r} \doublehat{X}_{im}^\ast \widetilde{S}_{mn}^0 \sum_{\ell=0}^{N_{\mu}-1} \doublehat{V}_{\ell n}^\ast  A_{0\ell}\label{DLRA-scheme: coupled equations v0}\\
&+ \frac{\Delta t}{1 + \sigma \Delta t} \sum_{i=1}^{N_x} \frac{1}{B_j^0} D_{ji}^{xx} B_i^0 \sum_{m,n=1}^{4r} \doublehat{X}_{im}^\ast \widetilde{S}_{mn}^0 \sum_{\ell=0}^{N_{\mu}-1} \doublehat{V}_{\ell n}^\ast  |A|_{0\ell} +\frac{\sqrt{2} \sigma \Delta t}{1 + \sigma \Delta t} \frac{B_j^1}{B_j^0} ,\nonumber\\
\frac{B_j^1}{B_j^0} =& \ 1 + \sigma \Delta t \frac{B_j^1}{B_j^0} \left(\sqrt{2} v_{j0}^1 - 2 \right).\label{DLRA-scheme: coupled equations B}
\end{align}
From \eqref{DLRA-scheme: coupled equations v0} and \eqref{DLRA-scheme: coupled equations B} we can then retrieve $\mathbf{v}_0^1 = \left(v_{j0}^1\right) $ and $ \mathbf{B}^1 = \left(B_j^1\right)$. We use the latter to divide out the factor $\frac{B_j^1}{B_j^0}$ in the low-rank representation of $\frac{B_j^1}{B_j^0} v_{jk}^\ast = \sum_{m,n=1}^{4r} \doublehat{X}_{jm}^\ast \doublehat{S}_{mn}^\ast \doublehat{V}_{kn}^\ast$. We perform the transformation step
\begin{align}\label{DLRA-scheme: correction step K}
    K_{jp}^{\ast,\text{trans}} = \frac{B_j^0}{B_j^1} K_{jp}^\ast.
\end{align}
From a QR-decomposition we obtain $\doublehat{\mathbf{X}}^{\ast, \text{trans}} \doublehat{\mathbf{S}}^{\ast,\text{trans}} = \text{qr}\left( \mathbf{K}^{\ast,\text{trans}}\right)$. Then we perform an additional basis augmentation step according to
\begin{align}\label{DLRA-scheme: basis augmentation mass conservation}
\doublehat{\mathbf{X}}^1 = \text{qr}\left( \left[\mathbf{v}_{0}^1, \doublehat{\mathbf{X}}^{\ast, \text{trans}} \right] \right), \quad \doublehat{\mathbf{V}}^1 = \text{qr}\left( \left[\mathbf{e}_1, \doublehat{\mathbf{V}}^\ast \right] \right),
\end{align}
where we add $\mathbf{v}_0^1$ to the updated spatial low-rank basis since the mass of the system is given by the zeroth order moment $\mathbf{v}_0^1$. In the directional basis, we add $\mathbf{e}_1 \in \mathbb{R}^{N_\mu}$, denoting the first unit vector in $\mathbb{R}^{N_\mu}$. This ensures mass conservation of the proposed low-rank scheme. Then, we have to adjust the coefficient matrix $\doublehat{\mathbf{S}}^{\ast,\text{trans}}$ correspondingly as 
\begin{align}\label{DLRA-scheme: S step correction mass conservation}
\doublehat{\mathbf{S}}^1 = \doublehat{\mathbf{X}}^{1,\top} \doublehat{\mathbf{X}}^{\ast, \text{trans}} \doublehat{\mathbf{S}}^{\ast, \text{trans}} \doublehat{\mathbf{V}}^\ast \left(\mathbf{I} - \mathbf{e}_1 \mathbf{e}_1^\top \right) \doublehat{\mathbf{V}}^1 + \doublehat{\mathbf{X}}^{1,\top} \mathbf{v}_0^1 \mathbf{e}_1^\top \doublehat{\mathbf{V}}^1 \in \mathbb{R}^{\left(4r+1 \right) \times \left(4r+1 \right)}.
\end{align}
In a last step, we truncate the augmented quantities $\doublehat{\mathbf{X}}^1, \doublehat{\mathbf{S}}^1$ and $\doublehat{\mathbf{V}}^1$ back to a new rank $r_1$, using the truncation strategy described in \cite{ceruti2022rank} or an adjusted truncation method inspired by \cite{einkemmerscalone2023} that ensures conservation of mass and will be given in Section \ref{sec5.4:Mass conservation}. Altogether, we obtain the updated low-rank factors $\mathbf{X}^1, \mathbf{S}^1$ and $\mathbf{V}^1$ such that $v_{jk}^1 = \sum_{m,n=1}^{r_1} X_{jm}^1 S_{mn}^1 V_{kn}^1$. The structure of the DLRA scheme is visualized in Figure \ref{fig:flowchart}. 

Note that this scheme is significantly different from the one presented in \cite{baumann2024} as the multiplicative structure leads to additional basis augmentations in \eqref{DLRA-scheme: X basis augmentation 3r} and \eqref{DLRA-scheme: V basis augmentation 3r} as well as a different solution of the coupled equations \eqref{DLRA-scheme: coupled equations v0} and \eqref{DLRA-scheme: coupled equations B}.
\end{subequations}

\begin{figure}[htp!]
    \centering
    \begin{tikzpicture}[node distance = 3.5cm,auto,font=\sffamily]
        \node[block](init){
        \textbf{input}
    \begin{varwidth}{\linewidth}\begin{itemize}
        \item internal energy at time $t_0$: $B_j^0$
        \item low-rank factors at time $t_0$:  $X^0_{jm},S^0_{mn}, V^0_{kn}$
        \item rank at time $t_0$: $r$
    \end{itemize}\end{varwidth}
        };
        \node[block, below = 0.35cm of init, fill=blue!20](KLstep){update bases according to \eqref{DLRA-scheme: K step} and \eqref{DLRA-scheme: L step}};
        \node[block, below = 0.75cm of KLstep,fill=red!20](augment1){augment bases with $X^0_{jm}, V^0_{kn}$};
        \node[block, below = 0.75cm of augment1,fill=red!20](augment1a){augment bases with $\frac{1}{B_j^0} \sum_{i=1}^{N_x} D_{ji}^x B_i^0 X_{im}^0$, $\frac{1}{B_j^0} \sum_{i=1}^{N_x} D_{ji}^{xx} B_i^0 X_{im}^0$\\
         and $\sum_{\ell=0}^{N_\mu} A_{k\ell} V_{\ell n}^0$, $\sum_{\ell=0}^{N_\mu} \left|A\right|_{k\ell} V_{\ell n}^0$
according to \eqref{DLRA-scheme: X basis augmentation 3r} and \eqref{DLRA-scheme: V basis augmentation 3r}};
        \node[block, below = 0.75cm of augment1a,fill=blue!20](SStep){update coefficient matrix according to \eqref{DLRA-scheme: S step}};
        \node[block, below = 0.75cm of SStep,fill=blue!20](u0Bupdate){update zeroth order moment and internal energy\\
        according to \eqref{DLRA-scheme: coupled equations v0} and \eqref{DLRA-scheme: coupled equations B}};
        \node[block, below = 0.75cm of u0Bupdate,fill=orange!20](correction){perform transformation step according to \eqref{DLRA-scheme: correction step K}};
        \node[block, below = 0.75cm of correction,fill=red!20](augment2){augment bases with $\mathbf{v}_0^1$ and $\mathbf{e}_1$ according to \eqref{DLRA-scheme: basis augmentation mass conservation}};
        \node[block, below = 0.75cm of augment2,fill=green!20](Scorrection){adjust coefficient matrix $\doublehat{S}_{mn}^{\ast,\text{trans}}$ according to \eqref{DLRA-scheme: S step correction mass conservation}};
        \node[block, below = 0.75cm of Scorrection,fill=yellow!20](truncate){truncate factors $\doublehat{ X}^1_{jm},\doublehat{S}^1_{mn},\doublehat{V}^1_{kn}$};
        \node[block, below = 0.35cm of truncate](out){
        \textbf{output}
    \begin{varwidth}{\linewidth}\begin{itemize}
        \item internal energy at time $t_1$: $B_j^1$
        \item low-rank factors at time $t_1$:  $X^1_{jm},S^1_{mn},V^1_{kn}$
        \item rank at time $t_1$: $r_1$
    \end{itemize}\end{varwidth}
        };
        \path[line] (init) -- node [near end] {} (KLstep);
        \path[line] (KLstep) -- node [near end] {} (augment1);
        \path[line] (augment1) -- node [near end] {} (augment1a);
        \path[line] (augment1a) -- node [near end] {} (SStep);
        \path[line] ([yshift=-0.1cm]init.east) -- ([xshift=0.35cm,yshift=-0.1cm] init.east) |- (augment1.east);
        \path[line] ([yshift=0.1cm]init.east) -- ([xshift=1.2cm,yshift=0.1cm] init.east) |- (augment1a.east);
        \path[line] (augment1a.west) -- ([xshift=-0.35cm] augment1a.west) |- (augment2.west);
        \path[line] (SStep) -- (u0Bupdate);
        \path[line] (u0Bupdate) -- node [near end] {} (correction);
        \path[line] (u0Bupdate.east) -- ([xshift=0.35cm] u0Bupdate.east) |- (augment2.east);
        \path[line] (correction) -- node [near end] {} (augment2);
        \path[line] (augment2) -- (Scorrection);
        \path[line] (Scorrection) -- node [near start] {} (truncate);
        \path[line] (truncate) -- node [near start] {} (out);
        \node[below right=0.0cm and 0.1cm of KLstep.south](test){$K^{\ast}_{jp}, L^{\ast}_{kp}$};
        \node[below right=-0.05 cm and 0.1cm of augment1.south](test){$\widehat{X}^{\ast}_{jm},\widehat{V}^\ast_{kn}$};
        \node[below right=-0.07cm and 0.1cm of augment1a.south](test){$\doublehat{X}^{\ast}_{jm},\doublehat{V}^\ast_{kn}$};
        \node[below right=-0.05cm and 0.1cm of SStep.south](test){$\doublehat{S}^{\ast}_{mn}$};
        \node[right=0.0cm and 0.6cm of correction.east](test){$B^{1}_{j}, v^{1}_{j0}$};
        \node[below right=0.0cm and 0.1cm of correction.south](test){$X_{jm}^{\ast, \text{trans}}, S_{mn}^{\ast, \text{trans}} $};
        \node[below right=-0.07cm and 0.1cm of augment2.south](test){$\doublehat{ X}_{jm}^1, \doublehat{ V}_{kn}^1$ };
        \node[below right=-0.07cm and 0.1cm of Scorrection.south](test){$\doublehat{ S}_{mn}^1$};
    \end{tikzpicture}
    \caption{Flowchart of the stable and conservative method \eqref{DLRA-scheme}.}
    \label{fig:flowchart}
\end{figure}

\subsection{Energy stability of the proposed low-rank scheme}

It can be shown that the proposed DLRA scheme preserves the energy stability of the full system. The rewriting of equations \eqref{eqs: Bg fully discrete} into \eqref{eqs: full reformulation for DLRA both} as well as the basis augmentations in \eqref{DLRA-scheme: X basis augmentation 3r} and \eqref{DLRA-scheme: V basis augmentation 3r} differentiate the proposed DLRA method from existing schemes and are crucial for the proof. 

\begin{theorem}
Under the time step restriction $\Delta t \leq \Delta x$, the fully discrete DLRA scheme \eqref{DLRA-scheme} is energy stable, i.e. it holds $E^1 \leq E^0$.
\end{theorem}
\begin{proof}
We start with the internal energy $\mathbf{B}$ and multiply \eqref{DLRA-scheme: coupled equations B} with $B_j^0 B_j^1$. This leads to 
\begin{align*}
\left(B_j^1 \right)^2 =& \ B_j^0 B_j^1 + \sigma \Delta t \left(B_j^1\right)^2 \left(\sqrt{2}v_{j0}^1 - 2 \right).
\end{align*}
Analogously to the proof of Theorem \ref{Theorem: Energy stability full system}, we rewrite the product $B_j^0 B_j^1$ and sum over $j$, rendering the expression
\begin{align}\label{Theorem3: B}
\frac{1}{2} \sum_{j=1}^{N_x} \left(B_j^1\right)^2 =& \ \frac{1}{2} \sum_{j=1}^{N_x} \left(B_j^0 \right)^2 - \frac{1}{2} \sum_{j=1}^{N_x} \left( B_j^1 - B_j^0\right)^2 + \sigma \Delta t \sum_{j=1}^{N_x} \left(B_j^1\right)^2 \left(\sqrt{2}v_{j0}^1 - 2 \right).
\end{align}
In the next step, we multiply \eqref{DLRA-scheme: S step} with $\doublehat{X}_{\alpha q}^\ast \doublehat{V}_{\beta p}^\ast$ and sum over $q$ and $p$. We introduce the projection operators $P_{\alpha j}^{X^\ast} = \sum_{q=1}^{4r} \doublehat{X}_{\alpha q}^\ast \doublehat{X}_{jq}^\ast$ and $P_{k\beta}^{V^\ast} = \sum_{p=1}^{4r} \doublehat{V}_{kp}^\ast \doublehat{V}_{\beta p}^\ast$ as well as the notations $v_{\alpha\beta}^\ast := \sum_{p,q=1}^{4r} \doublehat{X}_{\alpha q}^\ast \doublehat{S}_{qp}^\ast \doublehat{V}_{\beta p}^\ast$ and $v_{\alpha\beta}^0 := \sum_{p,q=1}^{4r} \doublehat{X}_{\alpha q}^\ast \widetilde{S}_{qp}^0 \doublehat{V}_{\beta p}^\ast$. We obtain
\begin{align}\label{Theorem3 proof: v ast}
v_{\alpha \beta}^\ast =& \ \frac{1}{1 + \sigma \Delta t} v_{\alpha \beta}^0 - \frac{\Delta t}{1 + \sigma \Delta t} \sum_{i,j=1}^{N_x} \sum_{k,\ell=0}^{N_\mu-1} P_{\alpha j}^{X^\ast} \frac{1}{B_j^0} D_{ji}^x B_i^0  v_{i\ell}^0 A_{k\ell} P_{k\beta}^{V^\ast}\\
&+ \frac{\Delta t}{1 + \sigma \Delta t} \sum_{i,j=1}^{N_x} \sum_{k,\ell=0}^{N_\mu-1} P_{\alpha j}^{X^\ast} \frac{1}{B_j^0}  D_{ji}^{xx} B_i^0   v_{i\ell}^0 |A|_{k\ell} P_{k\beta}^{V^\ast}\nonumber.
\end{align}
Further, we denote $\doublehat{v}_{\alpha\beta}^1 := \sum_{p,q=1}^{4r+1} \doublehat{X}_{\alpha q}^1 \doublehat{S}_{qp}^1 \doublehat{V}_{\beta p}^1$. From equation \eqref{DLRA-scheme: S step correction mass conservation}, we can derive the equation
\begin{align*}
\frac{B_\alpha^1}{B_\alpha^0} \doublehat{v}_{\alpha\beta}^1 = v_{\alpha \beta}^\ast \left( 1-\delta_{\beta 0} \right) + \frac{B_\alpha^1}{B_\alpha^0} v_{\alpha 0}^1 \delta_{\beta 0}.\nonumber
\end{align*}
We insert \eqref{Theorem3 proof: v ast} and \eqref{DLRA-scheme: coupled equations v0}, establishing the expression
\begin{align*}
\frac{B_\alpha^1}{B_\alpha^0} \doublehat{v}_{\alpha\beta}^1 \left(1+\sigma \Delta t \right) = \Bigg(&  v_{\alpha \beta}^0 - \Delta t \sum_{i,j=1}^{N_x} \sum_{k,\ell=0}^{N_\mu-1} P_{\alpha j}^{X^\ast} \frac{1}{B_j^0}  D_{ji}^x B_i^0  v_{i\ell}^0 A_{k\ell} P_{k\beta}^{V^\ast}\\
&+ \Delta t \sum_{i,j=1}^{N_x} \sum_{k,\ell=0}^{N_\mu-1} P_{\alpha j}^{X^\ast} \frac{1}{B_j^0}  D_{ji}^{xx} B_i^0   v_{i\ell}^0 |A|_{k\ell} P_{k\beta}^{V^\ast}\Bigg) \left(1-\delta_{\beta 0} \right)\\
+ \Bigg(& v_{\alpha \beta}^0 - \Delta t \sum_{i=1}^{N_x} \sum_{\ell=0}^{N_\mu-1} \frac{1}{B_\alpha^0} D_{\alpha i}^x B_i^0 v_{i\ell}^0 A_{0\ell}\\
&+ \Delta t \sum_{i=1}^{N_x} \sum_{\ell=0}^{N_\mu-1} \frac{1}{B_\alpha^0} D_{\alpha i}^{xx} B_i^0 v_{i\ell}^0 |A|_{0\ell} +  \sqrt{2} \sigma \Delta t \frac{B_\alpha^1}{B_\alpha^0}\Bigg)  \delta_{\beta 0}.
\end{align*}
We now use the fact that we have augmented the spatial basis according to \eqref{DLRA-scheme: X basis augmentation 3r} and \eqref{DLRA-scheme: V basis augmentation 3r}. This allows us to write any function $h_i^0 \in \text{span}\left( X_i^0 \right)$ and $\tilde{h}_\ell^0 \in \text{span}\left( V_\ell^0 \right)$ as
\begin{align*}
\sum_{i,j=1}^{N_x} P_{\alpha j}^{X^\ast} \frac{1}{B_j^0} D_{ji}^x B_i^0 h_i^0 = \frac{1}{B_\alpha^0} \sum_{i=1}^{N_x} D_{\alpha i}^x B_i^0 h_i^0 \hspace{0.25cm} \text{and} \hspace{0.25cm}\sum_{k,\ell=0}^{N_\mu-1} \tilde{h}_\ell^0 A_{k\ell} P_{k\beta}^{V^\ast} = \sum_{\ell=0}^{N_\mu-1} \tilde{h}_\ell^0 A_{\beta \ell},\\
\sum_{i,j=1}^{N_x} P_{\alpha j}^{X^\ast} \frac{1}{B_j^0} D_{ji}^{xx} B_i^0 h_i^0 = \frac{1}{B_\alpha^0} \sum_{i=1}^{N_x} D_{\alpha i}^{xx} B_i^0 h_i^0 \hspace{0.25cm} \text{and} \hspace{0.25cm}\sum_{k,\ell=0}^{N_\mu-1} \tilde{h}_\ell^0 \left|A\right|_{k\ell} P_{k\beta}^{V^\ast} = \sum_{\ell=0}^{N_\mu-1} \tilde{h}_\ell^0 \left|A\right|_{\beta \ell}.
\end{align*}
The basis augmentations as well as the properties of the projection operators enable us to obtain a representation of the form
\begin{align*}
\frac{B_\alpha^1}{B_\alpha^0} \doublehat{v}_{\alpha\beta}^1 \left(1+\sigma \Delta t \right) = F \left(1-\delta_{\beta 0} \right) + F  \delta_{\beta 0} +  \sqrt{2} \sigma \Delta t \frac{B_\alpha^1}{B_\alpha^0} \delta_{\beta 0}
\end{align*}
with
\begin{align*}
F = v_{\alpha \beta}^0 - \Delta t \sum_{i=1}^{N_x} \sum_{\ell=0}^{N_\mu-1} \frac{1}{B_\alpha^0} D_{\alpha i}^x B_i^0 v_{i\ell}^0 A_{\beta\ell} + \Delta t \sum_{i=1}^{N_x} \sum_{\ell=0}^{N_\mu-1} \frac{1}{B_\alpha^0} D_{\alpha i}^{xx} B_i^0 v_{i\ell}^0 |A|_{\beta\ell}.
\end{align*}
We use that on the right-hand side $F\delta_{\beta 0}$ cancels out. This yields the equation
\begin{align*}
\frac{B_\alpha^1}{B_\alpha^0} \doublehat{v}_{\alpha\beta}^1 \left(1+\sigma \Delta t \right) =& \ v_{\alpha \beta}^0 - \Delta t \sum_{i=1}^{N_x} \sum_{\ell=0}^{N_\mu-1} \frac{1}{B_\alpha^0} D_{\alpha i}^x B_i^0 v_{i\ell}^0 A_{\beta\ell} + \Delta t \sum_{i=1}^{N_x} \sum_{\ell=0}^{N_\mu-1} \frac{1}{B_\alpha^0} D_{\alpha i}^{xx} B_i^0 v_{i\ell}^0 |A|_{\beta\ell} + \sqrt{2} \sigma \Delta t \frac{B_\alpha^1}{B_\alpha^0} \delta_{\beta 0}.
\end{align*}
In the next step, we multiply with $B_\alpha^1 B_\alpha^0 \doublehat{v}_{\alpha\beta}^1$, sum over $\alpha$ and $\beta$, rearrange, and introduce the notation $\doublehat{u}_{\alpha\beta}^1 = B_\alpha^1 \doublehat{v}_{\alpha\beta}^1$. This leads to
\begin{align*}
\sum_{\alpha=1}^{N_x} \sum_{\beta=0}^{N_\mu-1} \left(\doublehat{u}_{\alpha\beta}^1\right)^2 =& \ \sum_{\alpha=1}^{N_x} \sum_{\beta=0}^{N_\mu-1} u_{\alpha \beta}^0 \doublehat{u}_{\alpha\beta}^1- \Delta t  \sum_{i,\alpha=1}^{N_x} \sum_{\ell,\beta=0}^{N_\mu-1}  \doublehat{u}_{\alpha\beta}^1 D_{\alpha i}^x u_{i\ell}^0 A_{\beta \ell}\\
&+ \Delta t \sum_{i,\alpha=1}^{N_x} \sum_{\ell,\beta=0}^{N_\mu-1}  \doublehat{u}_{\alpha\beta}^1 D_{\alpha i}^{xx} v_{i\ell}^0 |A|_{\beta \ell} +  \sigma \Delta t \sum_{\alpha=1}^{N_x} \sum_{\beta=0}^{N_\mu-1} \left(B_\alpha^1 \right)^2 \doublehat{v}_{\alpha\beta}^1 \left(\sqrt{2} \delta_{\beta 0} - \doublehat{v}_{\alpha\beta}^1 \right).
\end{align*}
Inserting the relation 
\begin{align*}
\sum_{\alpha=1}^{N_x} \sum_{\beta=0}^{N_\mu-1} u_{\alpha\beta}^0 \doublehat{u}_{\alpha\beta}^1 = \sum_{\alpha=1}^{N_x} \sum_{\beta=0}^{N_\mu-1} \left( \frac{1}{2} \left( \doublehat{u}_{\alpha\beta}^1 \right)^2 + \frac{1}{2} \left(u_{\alpha\beta}^0 \right)^2 - \frac{1}{2} \left( \doublehat{u}_{jk}^1 - u_{jk}^0 \right)^2 \right)
\end{align*}
establishes the expression
\begin{align*}
\frac{1}{2} \sum_{\alpha=1}^{N_x} \sum_{\beta=0}^{N_\mu-1} \left(\doublehat{u}_{\alpha\beta}^1\right)^2 =& \ \frac{1}{2} \sum_{\alpha=1}^{N_x} \sum_{\beta=0}^{N_\mu-1} \left(u_{\alpha\beta}^0\right)^2 - \frac{1}{2} \sum_{\alpha=1}^{N_x} \sum_{\beta=0}^{N_\mu-1} \left( \doublehat{u}_{jk}^1 - u_{jk}^0 \right)^2 - \Delta t  \sum_{i,\alpha=1}^{N_x} \sum_{\ell,\beta=0}^{N_\mu-1}  \doublehat{u}_{\alpha\beta}^1 D_{\alpha i}^x u_{i\ell}^0 A_{\beta \ell}\\
&+ \Delta t \sum_{i,\alpha=1}^{N_x} \sum_{\ell,\beta=0}^{N_\mu-1}  \doublehat{u}_{\alpha\beta}^1 D_{\alpha i}^{xx}    v_{i\ell}^0 |A|_{\beta \ell} +  \sigma \Delta t \sum_{\alpha=1}^{N_x} \sum_{\beta=0}^{N_\mu-1} \left(B_\alpha^1 \right)^2 \doublehat{v}_{\alpha\beta}^1 \left(\sqrt{2} \delta_{\beta 0} - \doublehat{v}_{\alpha\beta}^1 \right).
\end{align*}
This is the same expression as in equation \eqref{Theorem2-equ} in the proof of Theorem \ref{Theorem: Energy stability full system}. We apply the same estimates as in the proof of Theorem \ref{Theorem: Energy stability full system} and add the resulting equation with equation \eqref{Theorem3: B}. Analogously to the proof of Theorem \ref{Theorem: Energy stability full system} and due to the fact that the truncation step does not alter the zeroth order moment, we obtain energy stability of the multiplicative DLRA scheme under the time step restriction $\Delta t \leq \Delta x$.
\end{proof}

\subsection{Mass conservation}\label{sec5.4:Mass conservation}

In addition, the DLRA scheme \eqref{DLRA-scheme} can be shown to be mass conservative when using a suitable truncation strategy. We follow the ideas in \cite{einkemmerscalone2023, guo2024, einkemmerkuschschotthoefer2023} and adjust the truncation step such that it conserves the zeroth order moment. Different from \cite{einkemmerscalone2023} and as explained in \cite{einkemmerkuschschotthoefer2023}, we do not need to adjust the \textit{L}-step equation due to the usage of the augmented BUG integrator from \cite{ceruti2022rank}. Starting from the augmented quantities $\doublehat{\mathbf{X}}^1$, $\doublehat{\mathbf{S}}^1$ and $\doublehat{\mathbf{V}}^1$, the conservative truncation strategy then works as follows:
\begin{enumerate}
\item We set $\doublehat{\mathbf{K}}^1 = \doublehat{\mathbf{X}}^1 \doublehat{\mathbf{S}}^1$ and split it into two parts $\doublehat{\mathbf{K}}^1 = [\doublehat{\mathbf{K}}^{1,\text{cons}}, \doublehat{\mathbf{K}}^{1,\text{rem}}]$, where $\doublehat{\mathbf{K}}^{1,\text{cons}}$ corresponds to the first and $\doublehat{\mathbf{K}}^{1,\text{rem}}$ to the remaining columns of $\doublehat{\mathbf{K}}^1$. Analogously, we split $\doublehat{\mathbf{V}}^1$ into $\doublehat{\mathbf{V}}^1 = [\doublehat{\mathbf{V}}^{1,\text{cons}}, \doublehat{\mathbf{V}}^{1,\text{rem}}]$, where $\doublehat{\mathbf{V}}^{1,\text{cons}}$ corresponds to the first and $\doublehat{\mathbf{V}}^{1,\text{rem}}$ to the remaining columns of $\doublehat{\mathbf{V}}^1$. 

\item We compute $\doublehat{\mathbf{X}}^{1,\text{cons}} = \doublehat{\mathbf{K}}^{1,\text{cons}}/ \Vert \doublehat{\mathbf{K}}^{1,\text{cons}}\Vert$ and $\doublehat{\mathbf{S}}^{1,\text{cons}} = \Vert \doublehat{\mathbf{K}}^{1,\text{cons}}\Vert$. 

\item We perform a QR-decomposition of $\doublehat{\mathbf{K}}^{1,\text{rem}}$ to get $\doublehat{\mathbf{K}}^{1,\text{rem}} = \doublehat{\mathbf{X}}^{1,\text{rem}} \doublehat{\mathbf{S}}^{1,\text{rem}}$. 

\item We compute the singular value decomposition of $\doublehat{\mathbf{S}}^{1,\text{rem}} = \doublehat{\mathbf{P}} \mathbf{\Sigma} \doublehat{\mathbf{Q}}^\top$ with $\mathbf{\Sigma} = \diag(\sigma_j)$. The new rank $r_1 \leq 4r$ is chosen such that for a prescribed tolerance parameter $\vartheta$ it holds
\begin{align*}
\left(\sum_{j=r_1+1}^{4r} \sigma_j^2\right)^{1/2} \leq \vartheta.
\end{align*}
We set $\mathbf{S}^{1,\text{rem}} \in \mathbb{R}^{r_1\times r_1}$ to be the matrix containing the $r_1$ largest singular values. For the update of the spatial and the directional basis we introduce the matrices $\doublehat{\mathbf{P}}^{\text{rem}} \in \mathbb{R}^{4r\times r_1}$ and $\doublehat{\mathbf{Q}}^{\text{rem}} \in \mathbb{R}^{4r\times r_1}$ containing the first $r_1$ columns of $\doublehat{\mathbf{P}}$ and $\doublehat{\mathbf{Q}}$, respectively, and set $\mathbf{X}^{1,\text{rem}} = \doublehat{\mathbf{X}}^{1,\text{rem}} \doublehat{\mathbf{P}}^{\text{rem}} \in \mathbb{R}^{N_x \times r_1}$ and $\mathbf{V}^{1,\text{rem}} = \doublehat{\mathbf{V}}^{1,\text{rem}} \doublehat{\mathbf{Q}}^{\text{rem}} \in \mathbb{R}^{N_\mu \times r_1}$.

\item We set $\widetilde{\mathbf{X}}^1 = [\doublehat{\mathbf{X}}^{1,\text{cons}}, \mathbf{X}^{1,\text{rem}}]$ and $\widetilde{\mathbf{V}}^1 = [\mathbf{e}_1, \mathbf{V}^{1,\text{rem}}]$ and perform a QR-decomposition to obtain $\widetilde{\mathbf{X}}^1 = \mathbf{X}^1 \mathbf{R}^1$ and $\widetilde{\mathbf{V}}^1 = \mathbf{V}^1 \mathbf{R}^2$, respectively.

\item We compute 
\begin{align*}
\mathbf{S}^1 = \mathbf{R}^1 \begin{bmatrix}
\doublehat{\mathbf{S}}^{1,\text{cons}} & 0 \\ 
0 & \mathbf{S}^{1,\text{rem}}
\end{bmatrix} \mathbf{R}^{2,\top}.
\end{align*}
\end{enumerate}
This leads to the updated solution $\mathbf{v}^1 = \mathbf X^{1}\mathbf{S}^1\mathbf{V}^{1,\top}$ after one time step at time $t_1 = t_0+ \Delta t$.

In order to show local mass conservation for the proposed DLRA scheme, we translate the macroscopic quantities given in Definition \ref{Def1: Macroscopic quantities} to the fully discretized setting.

\begin{definition}[Fully discretized macroscopic quantities]
The \textit{mass} of the fully discretized multiplicative Su-Olson problem at time $t_0$ is defined as 
\begin{align*}
    \rho_j^0 := \sqrt{2} B_j^0 v_{j0}^0 + B_j^0. 
\end{align*}
The \textit{momentum} is given as 
\begin{align*}
u_j^0 := \sqrt{2} B_j^0 \sum_{\ell=0}^{N_\mu-1}v_{j\ell}^0 A_{0\ell}.
\end{align*}
For $t_1 = t_0 + \Delta t$ the definitions shall hold analogously.
\end{definition}
We can then show that the DLRA algorithm together with the conservative truncation strategy fulfills the following local conservation law.

\begin{theorem}
The DLRA scheme \eqref{DLRA-scheme} together with the conservative truncation strategy is locally mass conservative, i.e. it fulfills the local conservation law
\begin{align}\label{Theorem4: Local conservation law}
\frac{1}{\Delta t} \left(\sqrt{2} B_j^1 \Phi_j^1 + B_j^1 - \left( \sqrt{2} B_j^0 \Phi_j^0 + B_j^0\right) \right) = - \sqrt{2}\sum_{i=1}^{N_x} \sum_{\ell=0}^{N_\mu-1} D^x_{ji} B_i^0 v_{i\ell}^0 A_{0\ell}+ \sqrt{2} \sum_{i=1}^{N_x} \sum_{\ell=0}^{N_\mu-1} D^{xx}_{ji} B_i^0 v_{i\ell}^0 |A|_{0\ell}, 
\end{align}
where $\Phi_j^0 = \sum_{m,n=1}^r X_{jm}^0 S_{mn}^0 V_{0n}^0$, $\Phi_j^1 = \sum_{m,n=1}^{r_1} X_{jm}^1 S_{mn}^1 V_{0n}^1$ as well as $v_{jk}^0 = \sum_{m,n=1}^r X_{jm}^0 S_{mn}^0 V_{kn}^0$. This is a discretization of the continuous local conservation law given in \eqref{local conservation law}.
\end{theorem}
\begin{proof}
The conservative truncation strategy is designed to leave the zeroth order moment unchanged, i.e. it holds $\sum_{m,n=1}^{4r} \doublehat{X}_{jm}^1 \doublehat{S}_{mn}^1 \doublehat{V}_{0n}^1 = v_{j0}^1$. In addition, we know from the basis augmentation \eqref{DLRA-scheme: basis augmentation mass conservation} and the adjustment step \eqref{DLRA-scheme: S step correction mass conservation} that it holds $\sum_{m,n=1}^{4r} \doublehat{X}_{jm}^1 \doublehat{S}_{mn}^1 \doublehat{V}_{0m}^1 = \sum_{m,n=1}^{r_1} X_{jm}^1 S_{mn}^1 V_{0n}^1$. Combining both equalities, we obtain
\begin{align*}
\Phi_j^1 = \sum_{m,n=1}^{r_1} X_{jm}^1 S_{mn}^1 V_{0n}^1 = \sum_{m,n=1}^{4r} \doublehat{X}_{jm}^1 \doublehat{S}_{mn}^1 \doublehat{V}_{0n}^1 = v_{j0}^1.
\end{align*}
We insert this relation into the coupled equations \eqref{DLRA-scheme: coupled equations v0} and \eqref{DLRA-scheme: coupled equations B}. We multiply \eqref{DLRA-scheme: coupled equations v0} with $\sqrt{2}\left(1+\sigma \Delta t \right)$, rearrange it, and multiply both equations with $B_j^0$. This leads to
\begin{subequations}
\begin{align}
\label{eq: Conservation law g} \sqrt{2} B_j^1 \Phi_j^1 =& \ \sqrt{2} B_j^0 \Phi_j^0 - \sqrt{2} \Delta t \sum_{i=1}^{N_x} \sum_{\ell=0}^{N_\mu-1} D_{ji}^x B_i^0 \sum_{m,n=1}^{4r} \doublehat{X}_{im}^\ast \widetilde{S}_{mn}^0 \doublehat{V}_{\ell n}^\ast  A_{0\ell}\\
&+ \sqrt{2} \Delta t \sum_{i=1}^{N_x} \sum_{\ell=0}^{N_\mu-1} D_{ji}^{xx} B_i^0 \sum_{m,n=1}^{4r} \doublehat{X}_{im}^\ast \widetilde{S}_{mn}^0 \doublehat{V}_{\ell n}^\ast  |A|_{0\ell} + \sigma \Delta t B_j^1 \left(2 - \sqrt{2}\Phi_j^1 \right),\nonumber\\
B_j^1 =& \ B_j^0 + \sigma \Delta t B_j^1 \left(\sqrt{2} \Phi_j^1 - 2 \right).\label{eq: Conservation law B}
\end{align}
\end{subequations}
Due to the basis augmentations with $\mathbf{X}^0$ and $\mathbf{V}^0$ from the augmented BUG integrator it can be concluded that
\begin{align*}
\sum_{m,n=1}^{4r} \doublehat{X}_{im}^\ast \widetilde{S}_{mn}^0 \doublehat{V}_{\ell n}^\ast = \sum_{m,n=1}^r X_{im}^0 S_{mn}^0 V_{\ell n}^0 = v_{i\ell}^0. 
\end{align*}
We insert this relation into \eqref{eq: Conservation law g}, add equations \eqref{eq: Conservation law g} and \eqref{eq: Conservation law B}, and rearrange the result. This leads to the local conservation law \eqref{Theorem4: Local conservation law}, ensuring the local conservation of mass.
\end{proof}

\section{Numerical results}\label{sec6:Numerical results}

In this section, we compare the solution of the DLRA scheme \eqref{DLRA-scheme} to the solution of the full equations \eqref{eqs: full reformulation for DLRA both} to underline the efficiency and accuracy of the proposed method. We give different test examples in 1D and 2D that validate our theoretical results.

\subsection{1D plane source}\label{sec6.1:1D plane source}

We first examine the 1D plane source test case which is a common test example for thermal radiative transfer \cite{ganapol2001,ganapol2008,peng2020-2D}. We consider the spatial domain $D=[-10,10]$ and the directional domain $[-1,1]$. The initial condition is chosen to be the cutoff Gaussian
\begin{align*}
v\left(t=0,x \right) = \frac{1}{B^0} \max \left(10^{-4}, \frac{1}{\sqrt{2\pi\sigma_{\mathrm{IC}}^2}} \exp\left(-\frac{(x-1)^2}{2\sigma_{\mathrm{IC}}^2}\right)\right)
\end{align*}
with constant deviation $\sigma_{\mathrm{IC}}=0.03$. The traveling particles are initially centered around $x=1$ and move into all directions $\mu\in [-1,1]$. The initial internal energy is set to $B^0=1$ and the opacity to the constant value $\sigma=1$. For the low-rank computations we use an initial rank of $r=10$. The total mass $m^n$ at time $t_n$ is defined as $m^n = \Delta x \sum_{j=1}^{N_x} \left(\sqrt{2} B_j^n v_{j0}^n + B_j^n\right)$. As computational parameters we use $N_x=1000$ cells in the spatial domain and $N_\mu=500$ moments for the approximation in the directional domain. The time step size is determined by $\Delta t = \text{CFL} \cdot \Delta x$ with $\text{CFL}=0.99$, according to the corresponding CFL condition.

For the theoretical proof of energy stability of the multiplicative DLRA scheme \eqref{DLRA-scheme} the basis augmentations to rank $4r+1$ performed in \eqref{DLRA-scheme: X basis augmentation 3r} and \eqref{DLRA-scheme: V basis augmentation 3r} are essential. However, in practical implementations it can be observed that the standard basis augmentations to rank $2r+1$ provide similar solutions while being significantly faster.For this reason, we propose to leave out the basis augmentations presented in \eqref{DLRA-scheme: X basis augmentation 3r} and \eqref{DLRA-scheme: V basis augmentation 3r} in practical applications.

\begin{figure}[t]
    \centering
    \includegraphics[width = 1.0\linewidth]{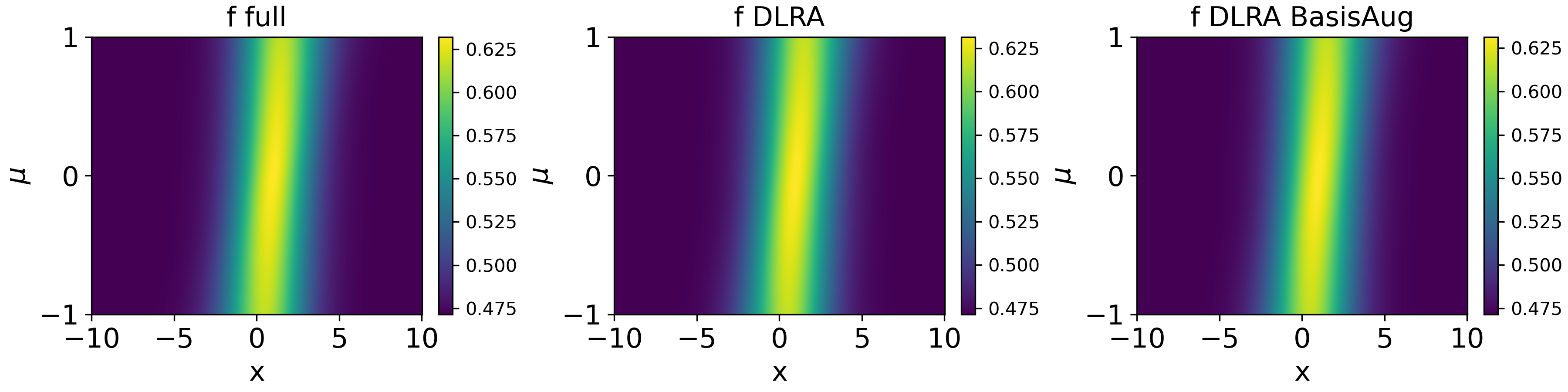}
   \caption{Numerical results for the solution $f(x,\mu)$ of the 1D plane source problem at time $t_{\text{end}}=8$ computed with the multiplicative full solver (left), the reduced multiplicative DLRA scheme (middle), and the augmented multiplicative DLRA scheme (right).}
    \label{fig:Planesource1}
\end{figure}

\begin{figure}[htb]
    \centering
    \includegraphics[width = 0.8\linewidth]{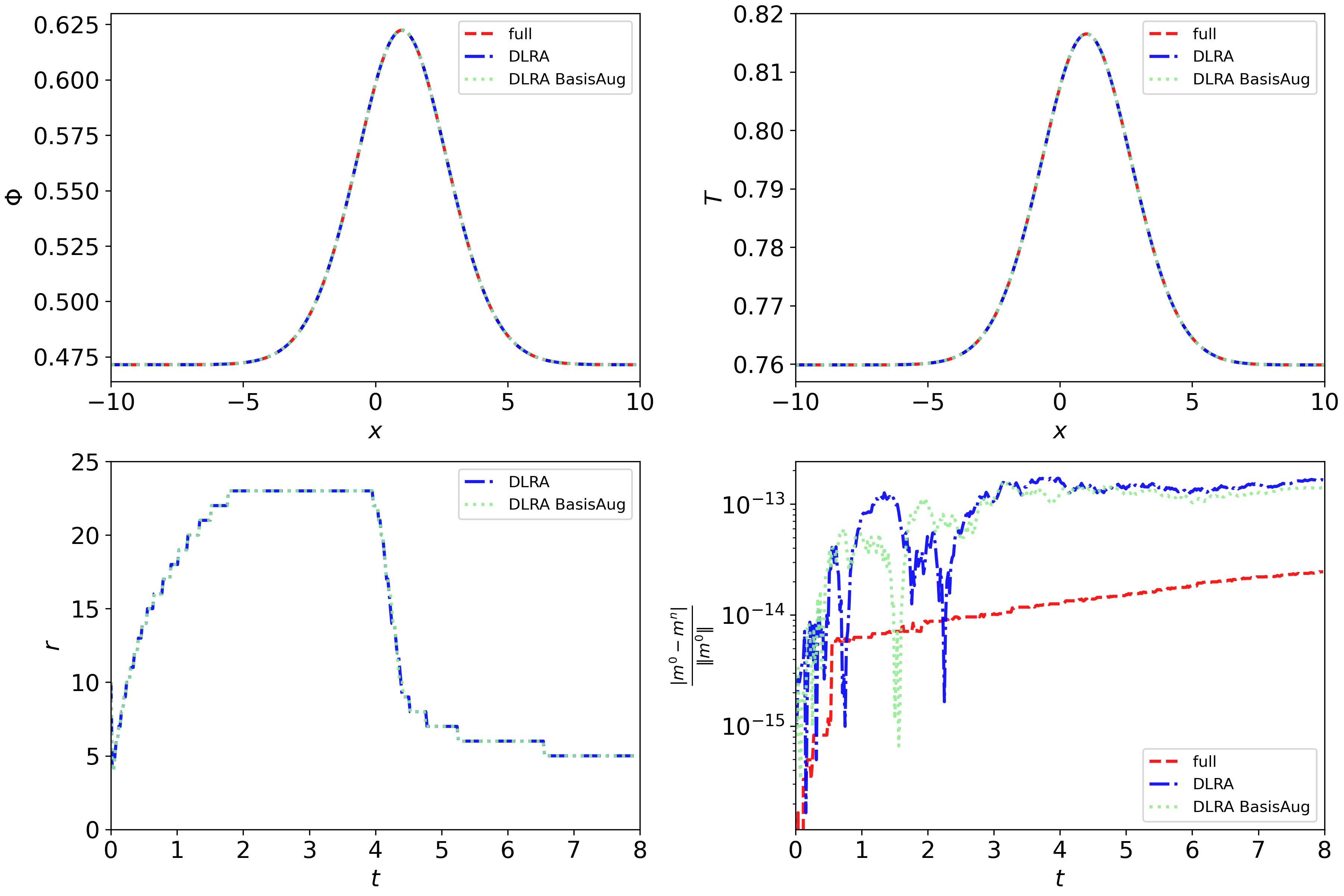}
   \caption{Top row: Numerical results for the scalar flux $\Phi$ (left) and temperature $T$ (right) of the 1D plane source problem at time $t_{\text{end}}=8$ computed with the multiplicative full solver, the reduced multiplicative DLRA scheme, and the augmented multiplicative DLRA scheme. Bottom row: Evolution of the rank in time for both multiplicative DLRA methods (left) and evolution of the relative mass error compared for all methods (right).}
    \label{fig:Planesource2}
\end{figure}

We compare the solution of the full system (full) with the solution of the reduced DLRA scheme with rank $2r+1$ (DLRA) and the solution of the augmented DLRA scheme with rank $4r+1$ (DLRA BasisAug). Figure \ref{fig:Planesource1} displays the solution $f(x,\mu)$ computed with the three different solvers. In Figure \ref{fig:Planesource2} the numerical results for the scalar flux $\Phi = \frac{1}{\sqrt{2}} \langle f \rangle_\mu$ and the dimensionless temperature $T = \sqrt[4]{B}$ are shown. All quantities are captured well by both DLRA schemes. In addition, Figure \ref{fig:Planesource2} shows the evolution of the rank $r$. For a chosen tolerance parameter of $\vartheta = 10^{-1} \Vert \mathbf{\Sigma} \Vert_F$ the rank of both DLRA schemes increases up to $r=23$ before it significantly reduces again. Note that the ranks of the reduced as well as the augmented DLRA scheme show a similar behavior as they are displayed after the corresponding truncation step. The relative mass error $\frac{\left|m^0-m^n \right|}{\left| m^0\right|}$ is of order $\mathcal{O}\left(10^{-13}\right)$, i.e. the DLRA scheme is mass conservative up to machine precision. These results confirm our theoretical findings. The memory demand reduce from $\mathcal{O} \left( N_x N_\mu \right)$ for the full solver to $\mathcal{O} \left( r \max \left(N_x, N_\mu \right) \right)$ for the DLRA schemes.

\subsection{1D external source}

In a second test example an external source term $Q(x)$ is added to the conservative form \eqref{eqs: Bg full} of the Su-Olson system, leading to the equations
\begin{align*}
\partial_t g (t,x,\mu) =& - \frac{\mu}{B(t,x)}\partial_x \left(B(t,x) g(t,x,\mu) \right) +\sigma \left(1-g(t,x,\mu) \right) - \frac{g(t,x,\mu)}{B(t,x)}\partial_t B(t,x) + \frac{Q(x)}{B(t,x)},\\
\label{eq: Bg full eq B}
\partial_t B(t,x) =& \ \sigma B(t,x)\left(\langle g(t,x,\mu) \rangle_{\mu}- 2 \right).
\end{align*}
The source term again generates radiation that moves through and interacts with the background material which in turn heats up and itself emits particles. The resulting travelling temperature wave is called a \textit{Marshak wave} \cite{marshak1958}. Again, we compare the solution of the full equations (full) with the solution of the reduced DLRA scheme with rank $2r+1$ (DLRA) and the solution of the augmented DLRA scheme with rank $4r+1$ (DLRA BasisAug). All schemes are adjusted to take the additional source term into account. For our numerical example we choose the function $Q(x) = \chi_{[-0.5,0.5]}(x)/a$ with $a= \frac{4 \sigma_{\mathrm{SB}}}{c}$ being the radiation constant. The initial internal energy is set to $B^0=50$. All other initial settings and computational parameters remain unchanged from the previous test example given in Section \ref{sec6.1:1D plane source}.

\begin{figure}[htb]
    \centering
    \includegraphics[width = 1.0\linewidth]{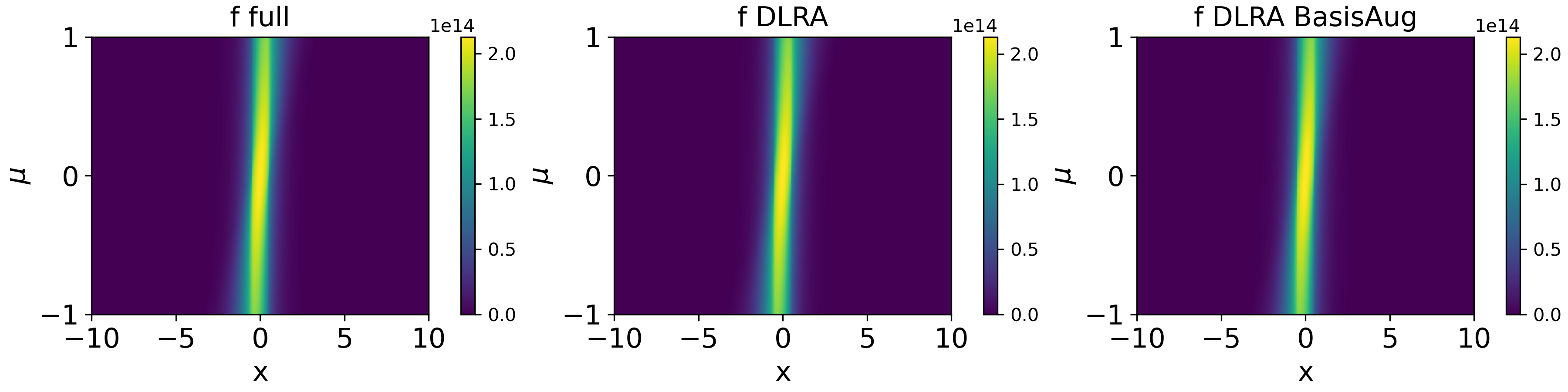}
   \caption{Numerical results for the solution $f(x,\mu)$ of the 1D external source problem at time $t_{\text{end}}=3.16$ computed with the  multiplicative full solver (left), the reduced multiplicative DLRA scheme (middle) and the augmented multiplicative DLRA scheme (right).}
    \label{fig:SuOlson1}
\end{figure}

\begin{figure}[htb]
    \centering
    \includegraphics[width = 1.0\linewidth]{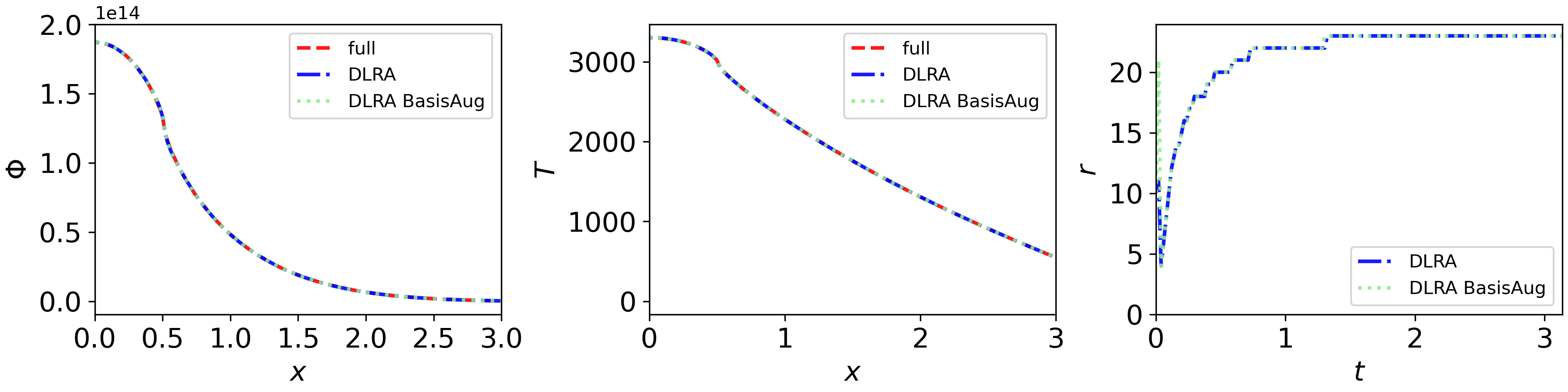}
   \caption{Numerical results for the scalar flux $\Phi$ (left) and temperature $T$ (middle) of the 1D external source problem at time $t_{\text{end}}=3.16$ computed with the multiplicative full solver, the reduced multiplicative DLRA scheme, and the augmented multiplicative DLRA scheme. Right: Evolution of the rank in time for both multiplicative DLRA methods.}
    \label{fig:SuOlson2}
\end{figure}

Figures \ref{fig:SuOlson1} and \ref{fig:SuOlson2} display the numerical results for the solution $f(x,\mu)$, for the scalar flux $\Phi = \frac{1}{\sqrt{2}} \langle f \rangle_\mu$, and for the dimensionless temperature $T=\sqrt[4]{B}$, computed with the full and both DLRA solvers. We again observe that both DLRA schemes capture the solution of the full system. The rank $r$ increases up to a value of $r=23$ for a chosen tolerance parameter of $\vartheta = 10^{-3} \Vert \mathbf{\Sigma} \Vert_F$. Due to the additional source term, there is no conservation of mass in this test example. 

\subsection{2D beam}

Computational benefits of the DLRA method can especially be observed in higher-dimensional settings. We extend the Su-Olson problem \eqref{eqs:Su-Olson} to a 2D spatial and a 2D angular setting. The corresponding equations is given by
\begin{align*}
\partial_t f(t,\mathbf{x},\mathbf{\Omega}) + \mathbf{\Omega} \cdot \nabla_\mathbf{x} f(t,\mathbf{x},\mathbf{\Omega}) &= \sigma(B(t,\mathbf{x})-f(t,\mathbf{x},\mathbf{\Omega})),\\
\partial_t B(t,\mathbf{x}) &= \sigma(\langle f(t,\mathbf{x},\cdot)\rangle_{\mathbf{\Omega}}-B(t,\mathbf{x})),
\end{align*}
where $\mathbf{x}=\left( x,y\right)$ and $\mathbf{\Omega} = \left( \Omega_x, \Omega_y, \Omega_z \right)$. For the numerical computations let $\mathbf{x} \in \left[-1,1 \right] \times \left[-1,1 \right]$ and $\mathbf{\Omega} \in \mathcal{S}^2$. The initial condition is chosen to be
\begin{align*}
f \left(t=0,\mathbf{x},\mathbf{\Omega} \right) = 10^6 &\cdot \frac{1}{2\pi\sigma_{\mathbf{x}}^2}\mathrm{exp}\left(-\frac{\Vert \mathbf{x} \Vert^2}{2\sigma_{\mathbf{x}}^2}\right)\cdot \frac{1}{2\pi\sigma_{\mathbf{\Omega}}^2}\mathrm{exp}\left(-\frac{(\Omega_x - \Omega^{\star} )^2 + (\Omega_z - \Omega^{\star} )^2}{2\sigma_{\mathbf{\Omega}}^2}\right),
\end{align*}
with $\Omega^{\star} = \frac1{\sqrt{2}}$, $\sigma_{\mathbf{x}} = \sigma_{\mathbf{\Omega}} = 0.1$. The initial internal energy is set to $B^0=1$ and the opacity to the constant value $\sigma=0.5$. The low-rank computations are performed with an initial rank of $r=100$. The total mass $m^n$ at time $t_n$ is defined as $m^n = \Delta x \Delta y \sum_{j=1}^{N_x \cdot N_y} \left(B_j^n v_{j0}^n + B_j^n\right)$. As computational parameters we use $N_x=N_y=500$ cells in each spatial dimension. For the 2D angular discretization we consider a spherical harmonics expansion with a polynomial degree of $N_{\mathbf{\Omega}}=29$, corresponding to $900$ expansion coefficients in angle. The time step size is determined by $\Delta t = \text{CFL} \cdot \Delta x$ with $\text{CFL}=0.7$.

\begin{figure}[t]
    \centering
    \includegraphics[width = 0.75\linewidth]{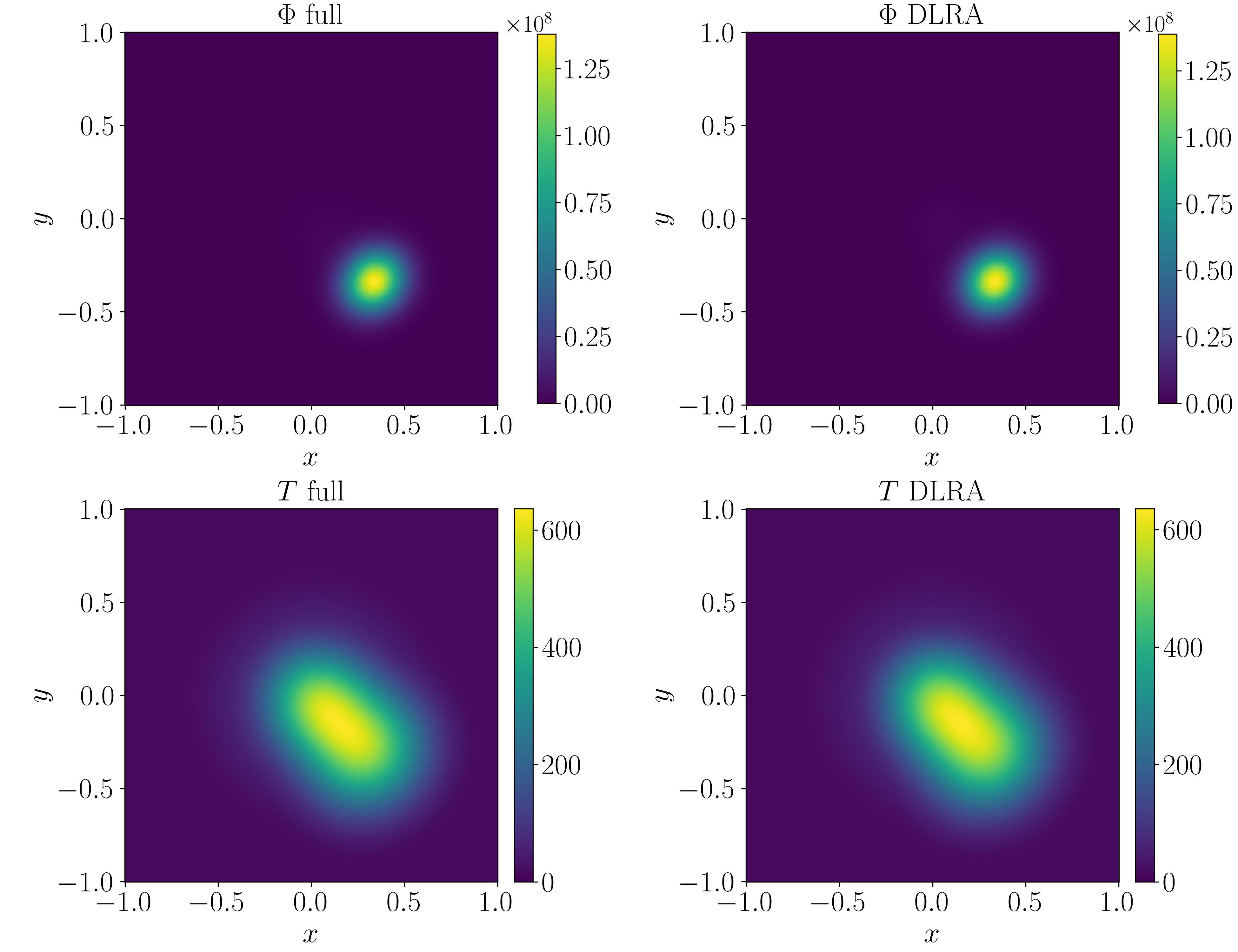}
   \caption{Numerical results for the scalar flux $\Phi$ (top row) and temperature $T$ (bottom row) of the 2D beam problem at time $t_{\text{end}}=0.5$ computed with the multiplicative full solver (left) and the reduced multiplicative DLRA scheme (right).}
    \label{fig:Beam1}
\end{figure}

\begin{figure}[t]
    \centering
    \includegraphics[width = 0.9\linewidth]{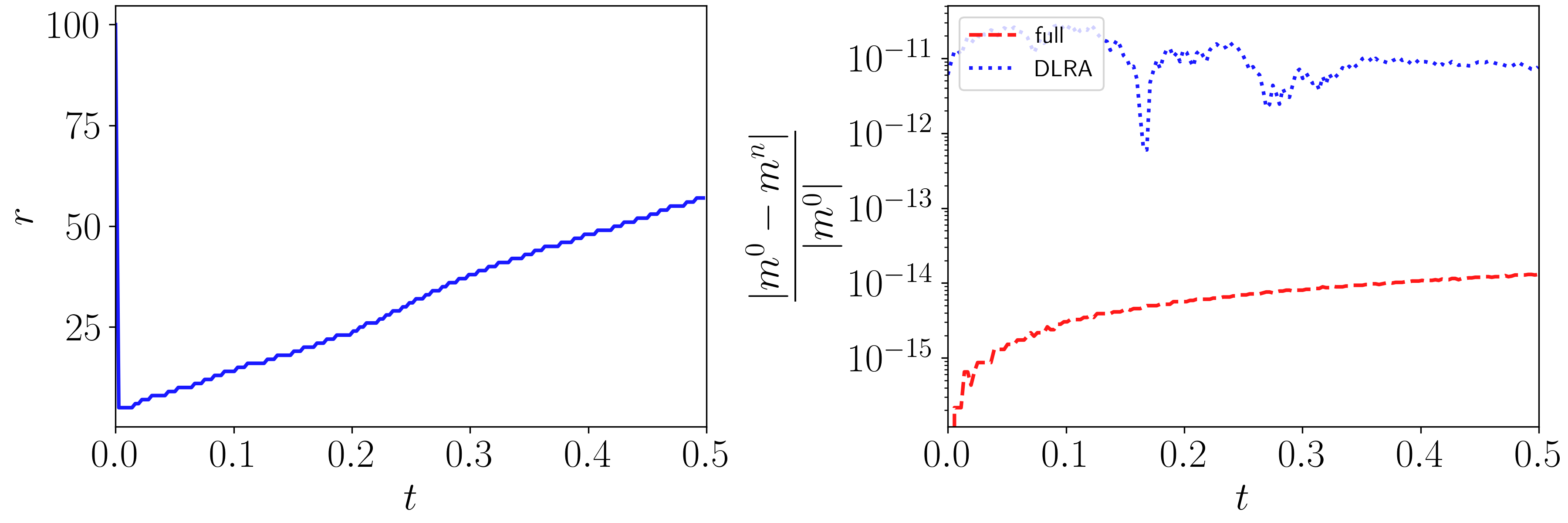}
   \caption{Evolution of the rank in time for the 2D beam problem for the reduced multiplicative DLRA method (left) and evolution of the relative mass error in time compared for both methods (right).}
    \label{fig:Beam2}
\end{figure}

In Figure \ref{fig:Beam1} the numerical results for the scalar flux $\Phi = \int_{\mathcal{S}^2} f(t,\mathbf{x}, \mathbf{\Omega}) \, \mathrm{d} \mathbf{\Omega}$ and the temperature $T = 4\pi \sqrt{2} \sqrt[4]{B}$ at time $t=0.5$, computed with full multiplicative solver and the reduced multiplicative DLRA scheme, are displayed. Note that the extension of the full as well as the DLRA scheme to 2D is straightforward. It can be observed that the solution of the full multiplicative solver is captured accurately by the multiplicative DLRA scheme. Figure \ref{fig:Beam2} presents the evolution of the rank $r$ of the reduced multiplicative DLRA scheme. For a chosen tolerance parameter of $\vartheta = 5 \cdot 10^{-4} \Vert \mathbf{\Sigma} \Vert_F$ the rank increases continuously over time but does not reach its allowed maximal value of $r_{\text{max}} = 100$. Both the full and the DLRA scheme are mass conservative as the relative mass errors are of order $\mathcal{O} \left(10^{-14} \right)$ and $\mathcal{O} \left(10^{-11} \right)$, respectively. For this setup, the computational benefit of the DLRA method is significant. The scheme is implemented in Julia v1.7 and performed on a MacBook Pro with M1 chip, resulting in a decrease of run time by a factor of approximately 35.5 from 61385 seconds to 1731 seconds. Due to the efficient implementation of the coupling of the equations in the multiplicative setting, the proposed multiplicative DLRA scheme is approximately 1.5 times faster than the non-multiplicative DLRA scheme presented in \cite{baumann2024}.

\section{Conclusion and outlook}\label{sec7:Conclusion}

We have presented a DLRA discretization for the multiplicative Su-Olson problem that is energy stable and mass conservative. To achieve both of these features, additional basis augmentations in the augmented BUG integrator combined with an adjusted truncation step are performed. This enables us to give a mathematically rigorous stability analysis. Numerical test examples confirm the theoretical properties and validate the accuracy and computational advantages of the DLRA scheme. However, the extension of the considered stability analysis from a linear to a non-linear problem, for example the isothermal Boltzmann-BGK equation treated in \cite{einkemmerhuying2021}, poses additional challenges as the general theoretical setting is significantly more difficult. For example, to our knowledge, the theoretical framework we have used here is not available in the non-linear case. Nevertheless, the analysis performed in this paper provides valuable insights into the choice of a suitable space discretization and stabilization when considering a multiplicative splitting of the distribution function, also for more complicated, potentially non-linear problems. In this sense, the theoretical considerations on the multiplicative splitting approach provided in this work can be extremely useful for the construction of stable and efficient DLRA schemes for more general equations, which we consider future work.

\section*{Acknowledgements}

Lena Baumann acknowledges support by the Würzburg Mathematics Center for Communication and Interaction (WMCCI) as well as the Stiftung der Deutschen Wirtschaft.

\newpage
\bibliographystyle{abbrv}
\bibliography{references}

\end{document}